\documentclass[a4paper,11pt]{amsart}
\def\JOURNAL{1}  %

\if\JOURNAL1
\usepackage{fullpage}
\fi
\usepackage{amsmath}
\usepackage{amsfonts}
\usepackage{amssymb}
\usepackage{listings}
\usepackage{mathtools}
\usepackage{hyperref}
\usepackage{amsthm}
\usepackage{verbatim}
\usepackage{bbm}

\usepackage{mathrsfs}%
\usepackage[title]{appendix}%
\usepackage{xcolor}%
\usepackage{textcomp}%
\usepackage{manyfoot}%
\usepackage{booktabs}%
\usepackage{algorithm}%
\usepackage{algorithmicx}%
\usepackage{algpseudocode}%
\usepackage{listings}%

\usepackage{todonotes}

\usepackage[T1]{fontenc}

\usepackage{array}
\usepackage{multirow}

\usepackage{caption}
\usepackage{subcaption}
\newcommand{\rset}{\mathbb{R}}
\newcommand{\rtset}{\mathbb{R}}
\newcommand{\tset}{\mathbb{R}}

\newcommand{\Bw}{\bar{\bar v}}
\newcommand{\Loss}{\mathcal{L}}
\newcommand{\PN}{\rho}
\newcommand{\HESS}{\mathrm{Hess}\,}

\def\XXint#1#2#3{{\setbox0=\hbox{$#1{#2#3}{\int}$ }
\vcenter{\hbox{$#2#3$ }}\kern-.6\wd0}}

\numberwithin{figure}{section}
\numberwithin{table}{section}
\numberwithin{equation}{section}
\newtheorem{theorem}{Theorem}[section]

\newtheorem{lemma}[theorem]{Lemma}

\newtheorem{preremark}[theorem]{Remark}

\title[Neural network approximation of Hamiltonian systems]{Convergence rates for random feature neural network approximation in molecular dynamics}

\if\JOURNAL1
\author[X. Huang]{Xin Huang}
\address{Institutionen f\"or Matematik, Kungl. Tekniska H\"ogskolan, 100 44 Stockholm, Sweden}
\email{xinhuang@kth.se}

\author[P. Plech\'a\v{c}] {Petr Plech\'a\v{c}}
\address{Department of Mathematical Sciences, University of Delaware, Newark, DE 19716, USA}
\email{plechac@udel.edu}

\author[M. Sandberg]{Mattias Sandberg}
\address{Institutionen f\"or Matematik, Kungl. Tekniska H\"ogskolan, 100 44 Stockholm, Sweden}
\email{msandb@kth.se}

\author[A. Szepessy]{Anders Szepessy}
\address{Institutionen f\"or Matematik, Kungl. Tekniska H\"ogskolan, 100 44 Stockholm, Sweden}
\email{szepessy@kth.se}
\date{}
\fi

\begin{document}

\if\JOURNAL1
\begin{abstract}
Random feature neural network approximations of the potential in Hamiltonian systems
yield  approximations of molecular dynamics correlation observables that 
has the expected error 
 $\mathcal{O}\big((K^{-1}+J^{-1/2})^{\frac{1}{2}}\big)$,
for  networks with $K$ nodes using $J$ data points,
provided the Hessians of the potential and the observables are bounded. 
The loss function is based on the least squares error of the potential and regularizations, with the data points sampled from the Gibbs density. The proof uses an elementary new derivation of the generalization error for random feature networks that does not apply the Rademacher or related complexities.
\end{abstract}


\fi

\maketitle

\section{Neural network approximations of Hamiltonian systems}
Approximation of 
potential energy surfaces 
by neural networks, trained by machine learning, has become a scientific success in
the sense that these approximate potentials can achieve  the
same computational accuracy as ab initio molecular dynamics but at a lower computational cost, see \cite{behler, weinan, review}.
Neural network based potentials therefore make it possible to study more complex systems that are 
out of reach for ab initio methods. Convergence rates for such neural network approximations have not been proved so far, which is the motivation for the present work.

We consider the Hamiltonian system
\begin{equation}\label{Ham}
\begin{split}
\dot x_t  &= p_t\,,\\
\dot p_t & = -\nabla V(x_t)\,,
\end{split}
\end{equation}
for $n$ particles in $\rset^3$, with the positions  $x:[0,\infty)\to \tset^{3n}$ and momenta $p:[0,\infty)\to\rset^{3n}$, 
where the particle masses are equal to one. 
Here $\dot x_t$ denotes the time
derivative of the position at time $t$, and similarly $\dot p_t$ for momenta. 
The function $V:\tset^{3n}\to\rset$ is a given potential and $\nabla V$ is the gradient of $V$, i.e., defining the force field.
The potential can be, for instance, a computational approximation 
to the ground state electron eigenvalue problem or 
the electron mean-field \cite{lebris,marx, xpma}.

We denote coordinates in the phase space
$\rset^{3n}\times\rset^{3n}$ by
\[
z_t:=
\left[
\begin{array}{c}
x_t\\ p_t\\
\end{array}
\right] \;\;\, \mbox{and define the vector field }\; 
f(z):=\left[\begin{array}{c}
p\\ -\nabla V(x)
\end{array} \right]\,,
\]
hence the Hamiltonian system \eqref{Ham} 
takes the form
\begin{equation}\label{eq:hamsysz}
\dot z_t = f(z_t)\,.
\end{equation}
The solution at time $t$ to \eqref{eq:hamsysz} with the initial condition $z_0$, at
$t=0$, is denoted $z_t(z_0)$.
We use the standard notation for the gradient of functions, i.e., 
$\nabla V$ for functions in the configuration space, 
and $\nabla_z A(x,p)\equiv (\nabla_x A,\nabla_p A)$ for functions in the phase space. 
We denote $\HESS V$ to be the Hessian of $V$, i.e., $(\HESS V)_{ij} := \partial^2_{x_ix_j} V$, and similarly $\HESS_z A$ is the Hessian of a function $A(z)$ in the phase space variables $z:= (x,p)$.

\medskip

The aim of this work is to derive error estimates 
for neural network approximations of molecular dynamics 
correlation observables 
\begin{equation}\label{obs}
\mathcal{C}_{AB}(t):=\int_{\rtset^{6n}} A\big(z_t(z_0)\big)\,B(z_0) \,\mu(z_0)\,\mathrm{d}z_0
\end{equation}
in the canonical ensemble, based on the Gibbs density 
\[
\mu(z):= \frac{1}{Z} e^{-\beta(|p|^2/2+V(x))}\,,\;\;\;
Z = \int_{\rtset^{6n}}  e^{-\beta(|p|^2/2+V(x))}\,\mathrm{d}x\,\mathrm{d}p\,,
\]
for given functions $A:\rtset^{6n}\to \rset$, $B:\rtset^{6n}\to \rset$, and a given inverse temperature $\beta>0$.

\medskip
We define the smooth cutoff functions $\chi_0$, $\chi_1\in\mathcal{C}^\infty(\rset^{3n})$, 
which for positive numbers $R_{-1}<R_0<R_1$, and
$j=0,1$ satisfy
\begin{equation}\label{smooth_cut_off_chi_j}
\chi_j(x)=\left\{
\begin{array}{ll}
1\,, & \mbox{ for }|x|<R_{j-1}\,,\\
0\,, & \mbox{ for }|x|> R_j\,.\\
\end{array}\right.
\end{equation}

Then we assume that the potential $V$ is confining in the sense that
it can be expressed in the form $V = \chi_0 V + v_e$ with 
an inner part $\chi_0 V$ and an outer part $v_e$ such that
the following assumptions are satisfied:
\begin{align}
    v_e(x) &:= (1-\chi_0(x))V(x)\,,\;\;\mbox{and} \;\;\int_{\rset^{3n}} e^{-\beta v_e(x)}\mathrm{d}x<\infty\,. \label{eq:ve} \\
   &\qquad \qquad \qquad \  \|\chi_0 V\|_\infty <\infty\,,  \label{eq:vi}
\end{align}%

\medskip\noindent
We assume that the confining part of the potential $v_e$ is known
and let 
\begin{equation}\label{eq:def_v1}
v:=\chi_1 V\,.
\end{equation}
We note that $\mathrm{supp}\chi_0 \subset \mathrm{supp}\chi_1$. We introduce  
a neural network approximation of the potential function $V$ which is constructed from
a random Fourier feature approximation of $v$. More precisely we define  
the approximation $\bar v$ which is a random Fourier 
feature neural network
\begin{equation}\label{eta_def}
\begin{split}
\bar v(x;\eta,\omega)&=\sum_{k=1}^K \eta_k  e^{\mathrm{i} \omega_k\cdot x}\,.
\end{split}
\end{equation}
The neural network approximation %
of the potential $V$ is then obtained
from the real part of $\bar v$ and the confining potential $v_e$ as 
\begin{equation}\label{eq:def_vr}
\begin{split}
  \Bar{V}&:= \chi_0 \Bar{v}+v_e\,,\\
  \bar v_r&:={\rm Re}(\bar V) = {\mathrm{Re}}(\chi_0\bar v)+ v_e = 
  \mathrm{Re}(\chi_0\bar v) +(1-\chi_0)V\,. %
\end{split}
\end{equation}
The neural network approximation 
$\bar v_r:\tset^{3n}\to\rset$
yields an approximating Hamiltonian system
\[
\begin{split}
\dot{\bar x}_t  &= \bar{p}_t\,,\\
\dot{\bar p}_t & = -\nabla \bar v_r(\bar x_t)\,,
\end{split}
\]
which we write as
\begin{equation}\label{eq:hamsysz_approx}
\dot{\bar z}_t = \bar f(\bar z_t)\,,
\end{equation}
using the notation
\[
\bar z_t:=
\left[
\begin{array}{c}
\bar x_t\\ \bar p_t\\
\end{array}
\right] \quad \mbox{ and }\quad 
\bar f(z):=\left[\begin{array}{c}
p\\ -\nabla \bar v_r(x)
\end{array} \right]\,.
\]
The purpose of the neural network approximation is 
to provide a computationally more efficient observable
\begin{equation}\label{C_ab_NN}
\Bar{\mathcal{C}}_{AB}(t):=\int_{\rtset^{6n}} A\big(\bar z_t(z_0)\big)\,B(z_0)\, \bar \mu(z_0)\,\mathrm{d}z_0\,,
\end{equation}
 compared to \eqref{obs}, where 
\[
\bar \mu(z):= \frac{1}{\bar Z}e^{-\beta(|p|^2/2+\bar v_r(x))}\,,\;\;\;
\bar Z = \int_{\rtset^{6n}}  e^{-\beta(|p|^2/2+\bar v_r(x))}\,\mathrm{d}z\,.
\]
The main question we study is: how accurate is the approximation \eqref{C_ab_NN} ?

\medskip
The central object of the approximating model potential $\Bar{v}_r$ is the  
random Fourier feature network function $\bar v$ which is based on $K$ 
nodes with the Fourier activation function,
where
$(\eta,\omega)\in \mathbb{C}^K\times\mathbb{R}^K$ are the network
parameters.
We denote $\|\eta\|_1 =\sum_{k=1}^K |\eta_k|$,
$\|\eta\|_2 = 
\big(\sum_{k=1}^K |\eta_k|^2\big)^{1/2}$ the usual
norms  on $\mathbb{C}^K$, and define the empirical loss function with $\alpha_1$, $\alpha_2>0$
\begin{equation}\label{fun:loss}
\ell(\eta,\omega) := 
        \frac{1}{J}\sum_{j=1}^J\big(
\alpha_1|v(x_j)-\bar{v}(x_j;\eta,\omega)|^2 + \alpha_2|\nabla v(x_j)-\nabla\bar{v}(x_j;\eta,\omega)|^2\big) 
\end{equation}
given the data set
\[
\{\big(x_j,V(x_j), \nabla V(x_j)\big)\ |\ j=1,\ldots,J\}\,,
\]
with the data points
$x_j$ independently sampled from the the Gibbs density
\[
\mu_x(x):=\frac{e^{-\beta V(x)}}{\int_{\rset^{3n}}
e^{-\beta V(x)}\,\mathrm{d}x}=\int_{\rset^{3n}}\mu(x,p)\,\mathrm{d}p\,.
\]

The standard learning problem for the neural network approximation 
$\bar{v}(\cdot;\eta,\omega)$ is then defined as the
statistical learning regression problem with the regularized 
empirical loss \eqref{fun:loss}
\begin{equation}\label{eq:learning}
\min_{(\eta,\omega)\in \mathbb{C}^K\times\mathbb{R}^K}
  \Big\{ \ell(\eta,\omega) 
  +\lambda_1\|\eta\|_2^2  +\lambda_2\|\eta\|_2^4 
+\lambda_3\max\big(\|\eta\|_1 -\tilde{C},0\big)\Big\}\,.
\end{equation}

\medskip

In this work we analyse a different approach which is motivated by randomized kernel 
methods. More specifically, the {\it random} Fourier feature %
network approximation 
differs from the minimization problem \eqref{eq:learning} 
by first drawing $K$ independent samples 
of $\omega_k\in \mathbb R^{3n}$ 
from a probability density $\PN:\mathbb R^{3n}\to [0,\infty)$, 
and then solving the optimization (regression) problem
\begin{equation}\label{opt_J}
\begin{split}
& \min_{\eta\in \mathbb C^K } %
\Big\{\frac{1}{J}\sum_{j=1}^J\big(
\alpha_1|v(x_j)-\bar{v}(x_j;\eta,\omega)|^2 + 
\alpha_2|\nabla v(x_j)-\nabla\bar{v}(x_j;\eta,\omega)|^2\big) 
\\&\quad  
+\lambda_1\sum_{k=1}^K|\eta_k|^2  +\lambda_2\big(\sum_{k=1}^K|\eta_k|^2\big)^2
+\lambda_3\max\big(\sum_{k=1}^K|\eta_k|-\tilde{C},0\big) \Big\}\,,
\end{split}
\end{equation} 
where $\tilde{C}$ is a sufficiently large constant 
chosen such that 
$\|\chi_1 V\|_{L^\infty}< \tilde{C}$.

How should the density $\rho$ for sampling
the frequencies in the random Fourier feature
be chosen?
The optimal sampling density for the frequencies
minimizes the variance of the learned regressor %
and is given by
\begin{equation}\label{rho_star_optimal}
    \rho_*(\omega):=
    \frac{|\widehat v(\omega)|\sqrt{1+|\omega|^2}}{\int_{\rset^{3n}}|\widehat v(\omega)|\sqrt{1+|\omega|^2}\,{\rm d}\omega}\,
\end{equation}
for the case with equal weight parameters $\alpha_1=\alpha_2$ in the loss function of \eqref{opt_J}. Here $\widehat v$ represents the Fourier transform of the inner part $v$ of the true potential function $V$, recalling that $v=\chi_1 V$.
To construct efficient random feature methods that provide 
independent frequencies from the optimal density 
is a challenge but approximate attempts exist, and
an algorithm has been developed in the adaptive random Fourier features method, see \cite{ajpma} and Section \ref{sec_numerics}.

We use the symbol
$\mathbb{E}$ 
for the expectation, in particular 
\[
\mathbb{E}_{{\omega}}[f] := \int_{\rset^{3n}} f(\omega)\PN(\omega) \mathrm{d}\omega 
\]
represents the expected value of $f:\mathbb{R}^{3n}\to\mathbb{C}$ under the distribution $\PN(\omega)\,\mathrm{d}\omega$,
and
\[
\mathbb{E}_{\{{x_j}\}}[h] := \int_{\rset^{3nJ}}h(x_1,\ldots,x_J)\mu_x(x_1){\rm d}x_1
\ldots\mu_x(x_J){\rm d}x_J%
\]
denotes the expectation of $h:\rset^{3nJ}\to\mathbb C$ with 
$x_1,\ldots,x_J$ independent, identically distributed (i.i.d.) samples
from the distribution $\mu_x(\mathrm{d}x)$.

The main result presented in this work is summarized in the following
theorem in which we derive the convergence rate for 
approximation of molecular dynamics observables
using neural network potentials.
\begin{theorem}\label{thm}
Given the potential $V=\chi_0 V+v_e$, and
$v=\chi_1 V$ with the splitting defined by
the smooth cut-off functions \eqref{eq:ve} and \eqref{eq:vi}, and denoting $\rho_*$ 
the optimal density \eqref{rho_star_optimal},
we assume that there exist a constant $C>0$ 
such that 
\begin{description}
\item[{\rm (i)}] 
the given external potential $v_e$ satisfies the bound 
\[
\begin{split}
\|\nabla v_e\|_{L^\infty} +\|\HESS v_e\|_{L^\infty}+\int_{\rset^{3n}}e^{-\beta v_e(x)}\,\mathrm{d}x& \leq C\,, 
\end{split}
\]
\item[{\rm (ii)}] the Fourier transform $\widehat v$ 
satisfies 
\[
\int_{\rset^{3n}}|\widehat v(\omega)|\,(1+|\omega|^2)^{\frac{1}{2}}\,\mathrm{d}\omega \leq C\,,
\]
\item[{\rm (iii)}] the observables $A$, $B$ and the sampling density $\PN$ satisfy
\begin{equation}\label{pn1}
\begin{split}
&\|A\|_{L^\infty} + \|\nabla_z A\|_{L^\infty} + \|\HESS_z(A)\|_{L^\infty} +\|B\|_{L^\infty} \leq C\,,\\
&\sup_{\omega\in\rset^{3n}} 
\frac{\rho_*(\omega)}{\PN(\omega)}\le C \,,\\
&\int_{\rset^{3n}}|\omega|^4\rho(\omega)\,\mathrm{d}\omega
\le C\,,
\end{split}
\end{equation}
\end{description}
then the  observable approximation \eqref{C_ab_NN},  
based on the neural network optimization \eqref{eta_def} 
and \eqref{opt_J} for
$\lambda_1=KJ^{-1/2}$, $\lambda_2=K^2J^{-1/2}$ and $\lambda_3=1$, has %
the expected error bound
\begin{equation}\label{C_b}
|\mathcal{C}_{AB}(t)-\mathbb E_{\{x_j\}}\big[\mathbb E_{\omega}[\bar{\mathcal{C}}_{AB}(t)]\big]|=
\mathcal O\Big(\big(K^{-1} + J^{-1/2}\big)^{1/2}\Big)\,,
\end{equation}
for a fixed $t<\infty$.
\end{theorem}

The convergence rate for the standard regression problem to approximate a given function by a neural network is well known \cite{barron, barron2, understand_ml}.  %
For instance, the work \cite{weinan_shallow} proves
that the 
generalization error, i.e., the mean squared error, with the probability $1-\delta$ 
has asymptotic behaviour 
$\mathcal O(K^{-1}+(J^{-1}\log( J/\delta))^{1/2})$ using $K$ nodes and $J$ data samples. 

The main new mathematical idea in the present work is to
construct and analyze the optimization problem \eqref{opt_J} combining  techniques  for
generalization error estimates with perturbation analysis of Hamiltonian systems. 
A new ingredient 
in our derivation  is that we prove in 
Theorem~\ref{thm:generalization} the generalization 
error estimate 
$\mathcal O( K^{-1}+J^{-1/2})$ %
without using the standard method based on Rademacher complexity or related concepts %
\cite{understand_ml,weinan_understand}. 
An alternative to the Rademacher complexity based on perturbation analysis is formulated in \cite{bous}. The analysis here also differs from \cite{bous}, which uses regularizations based on the Hilbert space norm induced by the reproducing kernel Hilbert space of random feature approximations. 
We use instead the regularization \eqref{opt_J} 
and the orthogonality provided by the random Fourier features, which in a certain sense separates the data dependence and the amplitudes $\eta$, although the amplitudes depend on the data.  

The next section gives an overview for the derivation of the main Theorem~\ref{thm}, and formulates an error estimate in Theorem \ref{thm:generalization} for the generalization error related to the optimization problem \eqref{opt_J}, using the random Fourier feature representation.
Section \ref{sec_numerics} presents numerical experiments on a computational model problem, with varied Fourier feature network size $K$ and dataset size $J$, demonstrating consistency with the theoretical error bounds. Section \ref{sec_proof} gives the formal proofs of Theorems \ref{thm} and \ref{thm:generalization}.

\section{Overview of main tools for proofs}
\subsection{The global/local error representation for Theorem~\ref{thm}} \label{global_local_error}
\if\JOURNAL1
\leavevmode  \\
\fi
To analyze the approximation error of a molecular dynamics observable
$g(z)$, we consider the function $u(z,t)$
defined by the transport equation
\begin{equation} \label{transport_eq_u}
\begin{split}
\partial_t u(z,t)+  f(z)\cdot \nabla_z  u(z,t) &=0\,,\quad t<T\,, \ z\in\rtset^{6n}\,,\\
 u(z, T) &= g(z)\,, \  \ z\in\rtset^{6n}\,,\\
\end{split}
\end{equation}
for a given observable $g:\rtset^{6n}\to\rset$. We have 
\[
 u(z,t) = g( z_T;  z_t=z)\,,
\]
where $z_s$ solves the evolution equation \eqref{eq:hamsysz} for $s\in[t,T]$.
Suppose $\bar z_s$ solves the approximate Hamiltonian dynamics \eqref{eq:hamsysz_approx} with neural network potential $\bar v_r$
then we obtain the following global/local error representation
\begin{equation}\label{g_est}
g(\bar z_T) -g(z_T) = \int_0^T \big(\nabla  V(\bar x_t)-\nabla \bar v_r( \bar x_t)\big)\cdot  \nabla_p  u(\bar z_t,t)\,\mathrm{d}t\,,
\end{equation}
since by using a telescoping cancellation we have
\begin{equation}\label{error_of_g}
    \begin{aligned}
    &\quad g(\bar{z}_T)-g(z_T)=u(\bar{z}_T,T)-u( z_0,0)\\&\stackrel{\{\text{$z_0=\bar{z}_0$}\}}{=}
    u(\bar{z}_T,T)-u(\bar{z}_0,0)
    =\int_0^T\mathrm{d}u(\bar{z}_t,t)\\
    &=\int_0^T \partial_t u(\bar{z}_t,t)+\dot{\bar{z}}_t\cdot\nabla_z u(\bar{z}_t,t)\ \mathrm{d}t
    =\int_0^T \partial_t u(\bar{z}_t,t)+\bar{f}(\bar{z}_t)\cdot\nabla_z u(\bar{z}_t,t)\ \mathrm{d}t\\
    &\stackrel{\text{by \eqref{transport_eq_u}}}{=}
    \int_0^T \big(-f(\bar{z}_t)+\bar{f}(\bar{z}_t)\big)\cdot \nabla_z u(\bar{z}_t,t)\,\mathrm{d}t\\
    &=\int_0^T \big(\nabla  V(\bar x_t)-\nabla \bar v_r( \bar x_t)\big)\cdot  \nabla_p  u(\bar z_t,t)\,\mathrm{d}t\,.
    \end{aligned}
\end{equation}
By considering the function $g$ as the observable $A$, the error representation \eqref{g_est} lays down the foundation for estimating the approximation error of the correlation observable $\Bar{\mathcal{C}}_{AB}$ in Theorem~\ref{thm}. The full proof of Theorem~\ref{thm} is provided in Section~\ref{proof_main_thm}.

\subsection{The generalization error for the optimization problem \eqref{opt_J}}
\if\JOURNAL1
\leavevmode \\
\fi
This subsection formulates the following estimate for the generalization error of the random Fourier feature optimization \eqref{eta_def} and \eqref{opt_J}. 
\begin{theorem}\label{thm:generalization} 
Given the potential $V=\chi_0 V+v_e$, and
$v=\chi_1 V$ with the splitting defined by
the smooth cut-off functions \eqref{eq:ve} and \eqref{eq:vi}, 
and denoting $\rho_*$ 
the optimal density \eqref{rho_star_optimal},
we assume that there exist a constant $C>0$ such that
\begin{description}
\item[{\rm (i)}] 
the given external potential $v_e$ satisfies $\int_{\rset^{3n}}e^{-\beta v_e(x)}\mathrm{d}x\leq C$,
\item[{\rm (ii)}] 
the Fourier transform $\widehat v$ of the inner potential $v$ satisfies %
\[\int_{\rset^{3n}}|\widehat v(\omega)|\,(1+|\omega|^2)^{\frac{1}{2}}\,\mathrm{d}\omega \leq C\,, %
\]
\item[{\rm (iii)}] 
the sampling density $\PN$ satisfies
\begin{equation}\label{pn}
\begin{split}
\sup_{\omega\in\rset^{3n}} 
\frac{\rho_*(\omega)}{\PN(\omega)}&\leq C \,,\\
\int_{\rset^{3n}}|\omega|^4\rho(\omega)\,\mathrm{d}\omega
&\le C\,,
\end{split}
\end{equation} %
\end{description}
then 
the optimization problem \eqref{eta_def} and \eqref{opt_J} 
has  the generalization error bound
\begin{equation}\label{gen_J}
\begin{split}
\mathbb E_{\{x_j\}}\Big[\mathbb E_\omega\big[ &
\int_{\rset^{3n}}\big(\alpha_1|v(x)-\bar{v}(x)|^2 + 
\alpha_2|\nabla v(x)-\nabla \bar{v}(x)|^2\big)\mu_x(x)\mathrm{d}x\big]\Big]\\
&=\mathcal O\big(K^{-1}+J^{-1/2}\big)\,,
\end{split}
\end{equation}
 and %
\[
\mathbb E_\omega[
\max(\sum_{k=1}^K|\eta_k|-\tilde C,0) ]
=\mathcal O(K^{-1})\,,
\] 
for $\lambda_1=K J^{-1/2}, \ \lambda_2=K^2J^{-1/2}$.
\end{theorem}

The complete proof of Theorem~\ref{thm:generalization} with the loss function $\mathbb{E}_\omega[|\bar{v}(x)-v(x)|^2+|\nabla \bar{v}(x)-\nabla v(x)|^2]$ is provided in Section~\ref{proof_thm_2_1}. Below we present a concise overview of the key ideas with the loss function $\mathbb{E}_\omega[|\bar{v}(x)-v(x)|^2]$.

We first introduce the notations $\mathbb E_x$ for the expectation with respect to the probability density $\mu_x$ and $\widehat{\mathbb E}_J$ for the empirical mean with respect to the data set  %
\begin{equation}\label{def_Ex_EJ}
\begin{split}
\mathbb E_x[|\bar{{v}}(x)-v(x)|^2]&:=\int_{\rset^{3n}}|\bar{{v}}(x)-v(x)|^2\,\mu_x(x)\,\mathrm{d}x\,,\\
\widehat{\mathbb E}_J[|\bar{{v}}(x)-v(x)|^2]&:=\frac{1}{J}\sum_{j=1}^J|\bar{{v}}(x_j)-v(x_j)|^2\,,\\
\end{split}
\end{equation}
and rewrite the generalization error as
\begin{equation}\label{generalization_error_v1}
\begin{split}
\mathbb E_{\{x_j\}}\Big[&\mathbb E_\omega\big[\mathbb E_x[|{\bar{v}}(x)-v(x)|^2\,\big|\,\{x_j \}]\big]\Big]
=\mathbb E_{\{x_j\}}\Big[\mathbb E_\omega\big[\widehat{\mathbb E}_J[|{\bar{v}}(x)-v(x)|^2\,\big|\,\{x_j \}]\big]\Big]\\
&+\mathbb E_{\{x_j\}}\Big[\mathbb E_\omega\big[(\mathbb E_x-\widehat{\mathbb E}_J)[|{\bar{v}}(x)-v(x)|^2\,\big|\,\{x_j \}]\big]\Big]\,.
\end{split}
\end{equation}
The first term of \eqref{generalization_error_v1} represents the training error, which has the bound 
\begin{equation}\label{training_error_estimate}
\begin{split}
&\quad\mathbb E_{\{x_j\}}\Big[\mathbb E_\omega\big[\widehat{\mathbb E}_J[|{\bar{v}}(x)-v(x)|^2\,\big|\,\{x_j \}]\\
&\qquad\qquad\quad\; +\lambda_1\sum_{k=1}^K|\eta_k|^2 + \lambda_2(\sum_{k=1}^K|\eta_k|^2)^2
+\lambda_3\max(\sum_{k=1}^K|\eta_k|-\tilde{C},0) 
\big]\Big]\\
&=\mathcal{O}\big(\frac{1}{K}+\frac{\lambda_1}{K}+ \frac{\lambda_2}{K^2}\big)\,.
\end{split}
\end{equation}
The derivation of \eqref{training_error_estimate} is given in Section~\ref{sub_sub_sec:training_error}. The primary challenge in estimating the remainder term of \eqref{generalization_error_v1} 
\[
\mathbb E_{\{x_j\}}\Big[\mathbb E_\omega\big[(\mathbb E_x-\widehat{\mathbb E}_J)[|{\bar{v}}(x)-v(x)|^2\,\big|\,\{x_j \}]\big]\Big]
\]
stems from the fact that the amplitude coefficients $\{\eta_k\}_{k=1}^K$ of the random Fourier feature representation $\Bar{v}(x_j)=\sum_{k=1}^K\eta_k \,e^{\mathrm{i}\omega_k\cdot x_j}$ are dependent on the specific training data set $\{(x_j,v(x_j)\}_{j=1}^J$, therefore
the difference of the expectation and the empirical mean $(\mathbb E_x -\widehat{\mathbb E}_J)[|{\bar{v}}(x)-v(x)|^2]$
does not reduce to a sum of independent variables with mean zero.
To address this issue,  we apply the Cauchy's inequality to split the error corresponding to the amplitude coefficients and the activation functions, and treat them separately.  A detailed analysis of the error estimate for the remainder term of \eqref{generalization_error_v1} %
is provided in Section~\ref{sub_sub_sec:testing_error_remainder}.

\subsubsection{The error estimate for the remainder term of \eqref{generalization_error_v1}}\label{sub_sub_sec:testing_error_remainder}
Using the notation for the trigonometric activation function
\[
\sigma(\omega\cdot x):= e^{\mathrm{i}\omega\cdot x}\,,
\]
we rewrite the remainder term of \eqref{generalization_error_v1} as
\begin{equation}\label{generalization_error_v1_square}
\begin{aligned}
    &  \mathbb{E}_{\{x_j \}}\Big[ \mathbb{E}_\omega\big[ (\mathbb{E}_x - \widehat{\mathbb{E}}_{J})[|{\bar{v}}(x)-v(x)|^2\,\big|\,\{x_j \}] \big]  \Big]\\
    & = \mathbb{E}_{\{x_j \}}\Big[ \mathbb{E}_\omega\big[ (\mathbb{E}_x - \widehat{\mathbb{E}}_{J})[ \sum_{k=1}^K \sum_{\ell=1}^K \eta_k^\ast\,\eta_\ell \,\sigma^\ast(\omega_k\cdot x)\,\sigma(\omega_\ell\cdot x)]\big]\Big] \\
    & \quad -2 \mathbb{E}_{\{x_j \}}\Big[ \mathbb{E}_\omega\big[ (\mathbb{E}_x - \widehat{\mathbb{E}}_{J})[ \sum_{k=1}^K  \mathrm{Re}\big(\eta_k \, \sigma(\omega_k\cdot x)\big)v(x)]\big]\Big]
    \\
    &\quad + \mathbb{E}_{\{x_j \}}\Big[ \mathbb{E}_\omega\big[ (\mathbb{E}_x - \widehat{\mathbb{E}}_{J})[ |v(x)|^2]\big]\Big]\,.
\end{aligned}
\end{equation}
For the first term of \eqref{generalization_error_v1_square}, we apply Cauchy's inequality to obtain
\begin{equation}\label{generalization_error_v1_Cauchy}
\begin{aligned}
      & \quad\ \textup{ } \big| \mathbb{E}_\omega\big[ \sum_{k=1}^K \sum_{\ell=1}^K \eta_k^\ast\, \eta_\ell\,\big( \widehat{\mathbb{E}}_J[ \sigma^\ast(\omega_k\cdot x)\sigma(\omega_\ell\cdot x) ] - \mathbb{E}_x[ \sigma^\ast(\omega_k\cdot x)\sigma(\omega_\ell\cdot x) ] \big) \big]    \big|\\
      & \leq \Big( \mathbb{E}_\omega\big[ \sum_{k=1}^K \sum_{\ell=1}^K |\eta_k\,\eta_\ell|^2 \big] \Big)^{\frac{1}{2}}\, \Big( \mathbb{E}_{\omega}\big[ \sum_{k=1}^K \sum_{\ell=1}^K \big| (\widehat{\mathbb{E}}_J - \mathbb{E}_x)[ \sigma^\ast(\omega_k\cdot x)\sigma(\omega_\ell\cdot x) ] \big|^2 \big] \Big)^{\frac{1}{2}}\\
      & \leq \frac{\lambda_2}{2} \mathbb{E}_\omega\big[ \big( \sum_{k=1}^K |\eta_k|^2 \big)^2 \big] + \frac{1}{2\lambda_2} \mathbb{E}_\omega\big[ 
\sum_{k=1}^K \sum_{\ell=1}^K \big| ( \widehat{\mathbb{E}}_J - \mathbb{E}_x )[ \sigma^\ast(\omega_k\cdot x)\sigma(\omega_\ell\cdot x) ] \big|^2 \big]\, ,
\end{aligned}
\end{equation}
where in the last inequality of \eqref{generalization_error_v1_Cauchy} we use that for any two real numbers $a$ and $b$ with a positive parameter $\gamma$ it holds that 
\[2|ab|\leq \frac{a^2}{\gamma}+\gamma b^2\,.\] 
By the training error estimate \eqref{training_error_estimate},
we have for the first regularization penalty term of \eqref{generalization_error_v1_Cauchy}
\[
    \frac{\lambda_2}{2} \mathbb{E}_{\{ x_j \}}\Big[ \mathbb{E}_\omega\big[ \big( \sum_{k=1}^K |\eta_k|^2 \big)^2 \big]  \Big] = \mathcal{O}\big( \frac{1}{K} + \frac{\lambda_1}{K} + \frac{\lambda_2}{K^2} \big)\,,
\]
and for the second term of \eqref{generalization_error_v1_Cauchy}, 
we obtain the expected error
\[
\mathbb{E}_{\{x_j\}}\Big[\mathbb{E}_\omega\big[ 
\sum_{k=1}^K \sum_{\ell=1}^K \big| ( \widehat{\mathbb{E}}_J - \mathbb{E}_x )[ \sigma^\ast(\omega_k\cdot x)\,\sigma(\omega_\ell\cdot x) ] \big|^2 \big]\Big] = \mathcal{O}\big(\frac{K^2}{J}\big)\,.
\]
The second term and the third term of \eqref{generalization_error_v1_square} can be treated similarly with the idea as above, which gives the generalization error estimate
\[
\begin{aligned}
&\quad\mathbb E_{\{x_j\}}\Big[\mathbb E_\omega\big[%
\int_{\rset^{3n}}%
|v(x)-\bar{v}(x)|^2 
\mu_x(x)\,\mathrm{d}x\big]\Big]\\
&=\,\mathcal O\big(\frac{1}{K}+\frac{\lambda_1}{K} + \frac{K}{\lambda_1 J} + \frac{\lambda_2}{K^2}+
\frac{K^2}{\lambda_2 J}\big)
=\mathcal O\big(\frac{1}{K} +\sqrt{\frac{1}{J}}\big)\,,\\
\end{aligned}
\]
for $\lambda_1=KJ^{-1/2}$ and $\lambda_2=K^2J^{-1/2}$. 

\subsubsection{The bound for the training error}\label{sub_sub_sec:training_error}

The derivation of the error estimate \eqref{training_error_estimate} uses 
the established method to apply Monte Carlo approximation of the inverse Fourier representation 
of the target function $v(x)$ as follows: consider
the Fourier transform %
\[
\widehat v(\omega) = \frac{1}{(2\pi)^{3n}}\int_{\tset^{3n}} v(x)\, e^{-\mathrm{i}\omega\cdot x}\,\mathrm{d}x\,,\ \omega\in\mathbb R^{3n}\,,
\]
and the inverse Fourier transform approximation
\begin{equation}\label{inverse_FT_v1}
v(x)=\int_{\mathbb R^{3n}} \widehat v(\omega)\,e^{\mathrm{i} \omega\cdot x}\,\mathrm{d}\omega
\simeq \sum_{k=1}^K \frac{\widehat v(\omega_k)\,e^{\mathrm{i} \omega_k\cdot x}}{K\, \PN(\omega_k)}
\,,\ x\in\tset^{3n}\,,
\end{equation}
with $K$ independent frequency samples $\{\omega_k\}_{k=1}^K$ drawn from the probability density $\PN$.
By defining the amplitude coefficients
\begin{equation} \label{eta_F}
\Hat{\eta}_k := \frac{\widehat v(\omega_k)}{K\, \PN(\omega_k)}\,,\ k=1,\dots,K,
\end{equation}
we obtain the corresponding approximation $\Bar{\Bar{v}}(x)$ of the potential function
\[ 
\Bar{\Bar{v}}(x) := \sum_{k=1}^K \Hat{\eta}_k \,e^{\mathrm{i} \omega_k\cdot x}\,,
\]
with the unbiased property
\[
\mathbb{E}_\omega[\bar{\bar{v}}(x)]=\sum_{k=1}^K \mathbb{E}_\omega\big[ \frac{\widehat v(\omega_k)}{ K\,\PN(\omega_k)} \,e^{\mathrm{i} \omega_k\cdot x} \big]=\frac{1}{K}\sum_{k=1}^K\int_{\mathbb R^{3n}} \widehat v(\omega_k)\,e^{\mathrm{i} \omega_k\cdot x}\,\mathrm{d}\omega_k=v(x)\,,
\]
and the variance 
\begin{equation}\label{squared_error_v_v1}
\mathbb{E}_\omega[|\bar{\bar{v}}(x)-v(x)|^2]=\frac{1}{K}\big( \int_{\mathbb{R}^{3n}}\frac{|\widehat v(\omega)|^2}{\PN(\omega)}\,\mathrm{d}\omega - |v(x)|^2 \big)=\mathcal{O}(K^{-1})\,.
\end{equation}
The variance \eqref{squared_error_v_v1} is minimized with the optimal probability density function
\[ \PN(\omega) = \frac{|\widehat v(\omega)|}{\| \widehat v(\omega)  \|_{L^1}}\,. \]
For the expected loss involving also the force fields $\mathbb{E}_\omega[|\bar{\bar{v}}(x)-v(x)|^2+|\bar{\bar{v}}'(x)-v'(x)|^2]$, we apply similarly the Fourier transform of $v(x)$ and $v'(x)$ to obtain the optimal probability density $\rho_\ast(\omega)=|\widehat v(\omega)|\sqrt{1+|\omega|^2}/\|\widehat v(\omega)\sqrt{1+|\omega|^2}\|_{L^1}$ which minimizes the variance of the approximation as given in \eqref{rho_star_optimal}.
Here for convenience, we use the shorthand notation $v'$ to represent the gradient of the function $v$.

The optimization \eqref{opt_J}, the bound \eqref{pn}, and the variance estimate \eqref{squared_error_v_v1} provide the following estimate for the training error
\begin{equation*}\label{training_error_v111}
\begin{split}
&\mathbb E_\omega\Big[\min_{\eta\in \mathbb C^K }\Big( %
\widehat{\mathbb E}_J[|{\bar{v}}(x)-v(x)|^2] +\lambda_1\sum_{k=1}^K|\eta_k|^2 + \lambda_2(\sum_{k=1}^K|\eta_k|^2)^2 
\\
& \qquad \qquad \qquad \qquad \qquad \qquad \quad \;
+\lambda_3\max(\sum_{k=1}^K|\eta_k|-\tilde{C},0) 
\Big)\Big]\\
&\le \mathbb E_\omega\big[\widehat{\mathbb E}_J[|\bar{\bar{v}}(x)-v(x)|^2] 
+\lambda_1\sum_{k=1}^K|\hat\eta_k|^2 + \lambda_2(\sum_{k=1}^K|\hat\eta_k|^2)^2
+\lambda_3\underbrace{\max(\sum_{k=1}^K|\hat\eta_k|-\tilde{C},0)}_{=0}\big]\\
&=\mathcal O(\frac{1}{K}+\frac{\lambda_1}{K}+ \frac{\lambda_2}{K^2})
\,.\\
\end{split}
\end{equation*}

\section{Numerical Experiments}\label{sec_numerics}
In this section, we conduct a systematic investigation on a model problem, examining the generalization error of the Fourier network representation. We evaluate its performance in reconstructing the potential field from a given training dataset, as well as in approximating correlation observables, and provide consistent numerical experiments with the theoretical estimates.

Section~\ref{computational_model}  outlines our computational model for the target potential function, $V(x)$. In Section~\ref{subsec_nume_general_error} we examine the behavior of the generalization error as the number of nodes $K$ in the Fourier network increases, corresponding to the error bound in Theorem~\ref{thm:generalization}. Subsequently, in Section~\ref{subsec_nume_corr_func}, we evaluate the approximation of the correlation observables by the trained Fourier network. Section~\ref{subsec_nume_one_beta} compares two sampling strategies for obtaining a comprehensive training data set. Lastly, Section~\ref{subsec:optimization_alg} provides more details of our implementation of the optimization algorithm for the training of the Fourier network.

\subsection{The computational model } \label{computational_model}
\if\JOURNAL1
\leavevmode \\
\fi
In order to test the prediction ability of the trained Fourier network and to visualize our results effectively, we devise a target potential function $V$ for ${x}\in\mathbb{R}^2$
\[
V({x})=\frac{1}{2}x_1^2+\frac{\alpha}{2}x_2^2+\gamma \sin{(x_1x_2)}\chi_{a}(x)\,,
\]
with $\alpha=\sqrt{2}$, $\gamma=2$, which is a variant of the potential surface leading to ergodic molecular dynamics in \cite{crossing_potential}. The inner part $v({x})$ of the target function $V({x})$ is given
by
\begin{equation}
    v({x})=V({x})\,\chi_c(x)\, ,
\end{equation}
where $\chi_a(x)$ and $\chi_c(x)$ are two smooth cutoff functions defined by
\[
\chi_{a}(x)=\left\{
\begin{array}{ll}
1\,, & \mbox{ for }|x|\le R_{a}\,,\\
e^{-\frac{(|x|-R_{a})^2}{2\delta_a}}\,, & \mbox{ for }|x|> R_{a}\,,\\
\end{array}\right.\quad
\chi_{c}(x)=\left\{
\begin{array}{ll}
1\,, & \mbox{ for }|x|\le R_{c}\,,\\
e^{-\frac{(|x|-R_{c})^2}{2\delta_c}}\,, & \mbox{ for }|x|> R_{c}\,,\\
\end{array}\right.
\]
imitating the smooth cutoff functions $\chi_j(x)$ in \eqref{smooth_cut_off_chi_j}. In practice we take the parameters $R_a=2$, $\delta_a=1$, $R_c=4$, $\delta_c=0.2$. 

In practise it is possible to generate $x$-data samples from the Gibbs density $\mu_x$ by the Metropolis method or numerical approximation of the overdamped Langevin dynamics 
\begin{equation}\label{langevin}
\mathrm{d}x_s=-\nabla V(x_s)\,\mathrm{d}s +\sqrt{\frac{2}{\beta}}\,\mathrm{d}W_s\,,
\end{equation}
where $W:[0,\infty)\to\rset^{3n}$ is the standard Wiener process,
using for instance the Leimkuhler-Matthews method \cite{Leimkuhler}, which at each time step and nuclei position $x$
requires the forces $\nabla V(x)$. Such nuclei positions at each time step will not be independent but after sufficiently
long time, exceeding the correlation length, the positions will be approximately independent due to the ergodic property of the dynamics.

\subsection{Generalization error of the trained Fourier network} \label{subsec_nume_general_error}
\if\JOURNAL1
\leavevmode \\
\fi
We aim to train a Fourier network representation 
\[\bar{v}(x; \omega, \eta)=\sum_{k=1}^K \eta_k \, e^{\mathrm{i} \omega_k\cdot x}\]
to approximate the unknown inner part of the potential function $v(x)$, based on the training data set
$\{\big(x_j,v(x_j), v'(x_j)\big)\}_{j=1}^J$
sampled with overdamped Langevin dynamics \eqref{langevin}
under the target potential function $V(x)$.

Utilizing $K$ independent and identically distributed samples of frequency parameters $\omega=\{\omega_k\}_{k=1}^K$, the random Fourier feature method trains the amplitude coefficients $\eta=\{\eta_k\}_{k=1}^K$ by solving the optimization problem \eqref{opt_J} with respect to $\eta_k\in\mathbb{C}$ for fixed frequency samples $\omega_k\in\mathbb{R}^2$, $k=1,2,\dots, K$.

On the other hand, the standard Fourier neural network treats both the frequencies $\omega$ and their corresponding amplitudes $\eta$ as trainable parameters, aiming to minimize the following regularized objective function $\mathcal{L}_{\mathrm{R}}(\omega,\eta)$ based on the loss function $\ell(\eta, \omega)$ as defined in \eqref{fun:loss}
\begin{equation}\label{regularized_loss}
\begin{aligned}
\mathcal{L}_{\mathrm{R}}(\omega,\eta)&:=
\frac{1}{J}\sum_{j=1}^J\big(
\alpha_1|v(x_j)-\bar{v}(x_j;\omega,\eta)|^2 + \alpha_2|v'(x_j)-\bar{v}'(x_j;\omega,\eta)|^2\big) 
\\&\quad  
+\lambda_1\sum_{k=1}^K|\eta_k|^2  +\lambda_2(\sum_{k=1}^K|\eta_k|^2)^2
+\lambda_3\max(\sum_{k=1}^K|\eta_k|-\tilde{C},0) 
\,,\\
\end{aligned}
\end{equation}
with respect to $\omega$ and $\eta$. Optimization of the frequencies $\omega$ can be achieved through the \emph{gradient descent} algorithm by evaluating the gradient of $\mathcal{L}_{\mathrm{R}}(\omega,\eta)$ with respect to $\omega$. Additionally, we implement the \emph{adaptive Metropolis} algorithm \cite{ajpma} to iteratively update the distribution of frequency parameters $\omega$, which aims to equidistribute the corresponding amplitudes.

Due to the presence of the last penalization term in \eqref{regularized_loss} which leads to a similarity to the Lasso regression \cite{Lasso, Lasso_paper}, we implement the proximal gradient method \cite{Prox_Grad_method} to find the optimal amplitude coefficients $\eta$ 
for a set of fixed frequencies $\omega$. Then the gradient $\nabla_{\omega\,} \mathcal{L}_{\mathrm{R}}$ is evaluated, facilitating the update on the frequency parameters $\omega$ with the 
standard gradient descent algorithm. 
Thus we iteratively update on $\eta$ and $\omega$ to achieve the optimization of $\mathcal{L}_{\mathrm{R}}$. For all the numerical tests in Section~\ref{sec_numerics}, we use the parameters $\alpha_1=1$, $\alpha_2=1$, $\lambda_1=10^{-2}$, $\lambda_2=10^{-3}$, $\lambda_3=10^{-2}$, and $\tilde{C}=100$. More details for our implementation of the optimization algorithms are summarized in Section ~\ref{subsec:optimization_alg}.

Corresponding to the estimate in Theorem~\ref{thm:generalization} for the generalization error, we evaluate the empirical testing loss 
\[
\begin{split}
\mathcal{L}_{\mathrm{emp}}&:=\frac{1}{Q}\sum_{q=1}^Q\Big(\frac{1}{J}\sum_{j=1}^J\big(
\alpha_1|v(\tilde{x}_j^{(q)})-\bar{v}(\tilde{x}_j^{(q)};\omega^{(q)},\eta^{(q)})|^2\\
&\qquad \qquad \qquad\quad\; \ + \alpha_2|v'(\tilde{x}_j^{(q)})-\bar{v}'(\tilde{x}_j^{(q)};\omega^{(q)},\eta^{(q)})|^2\big) \Big),
\end{split}
\]
based on $Q$ replicas of the training procedure using independent \emph{training data sets} \\
$\big\{\big(x_j^{(q)}, v(x_j^{(q)}),\nabla v(x_j^{(q)})\big)\big\}_{j=1}^J$ 
which leads to the corresponding optimized parameter sets $\omega^{(q)}$ and $\eta^{(q)}$ for $q=1,2,\dots,Q$. \\
For each replica an independent \emph{testing data set} $\big\{\big(\tilde{x}_j^{(q)}, v(\tilde{x}_j^{(q)}),\nabla v(\tilde{x}_j^{(q)})\big)\big\}_{j=1}^J$ is sampled under the same distribution as the training data set using the overdamped Langevin dynamics \eqref{langevin}, and is used to compute the empirical loss $\mathcal{L}_{\mathrm{emp}}$ serving as an approximation to the expected generalization error with respect to the data and with respect to the frequencies
\[
 \mathcal{L}_{\mathrm{emp}}\approx \mathbb E_{\{x_j\}}\Big[\mathbb E_\omega\big[%
\int_{\rset^{3n}}\big(\alpha_1|v(x)-\bar{v}(x)|^2 + \alpha_2|v'(x)-\bar{v}'(x)|^2\big)\mu_x(x)\,\mathrm{d}x\big]\Big]\,.
\]

\begin{figure}[ht!] 
  \centering
  \begin{subfigure}{0.51\textwidth}
    \includegraphics[width=\linewidth]{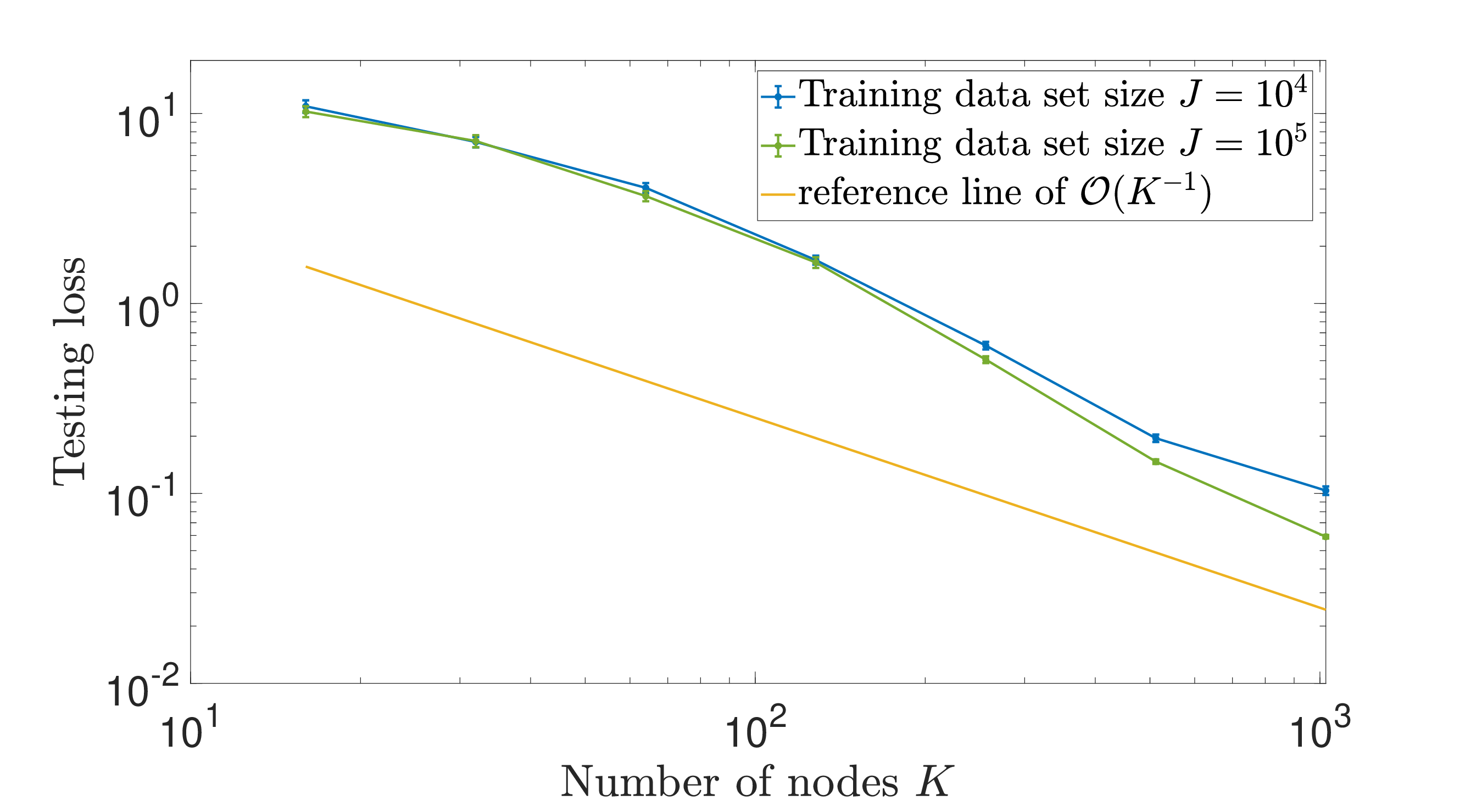}
    \caption{ Adaptive Metroplis for $\omega$ }
  \end{subfigure}
  \hspace{-0.04\textwidth} %
  \begin{subfigure}{0.51\textwidth}
    \includegraphics[width=\linewidth]{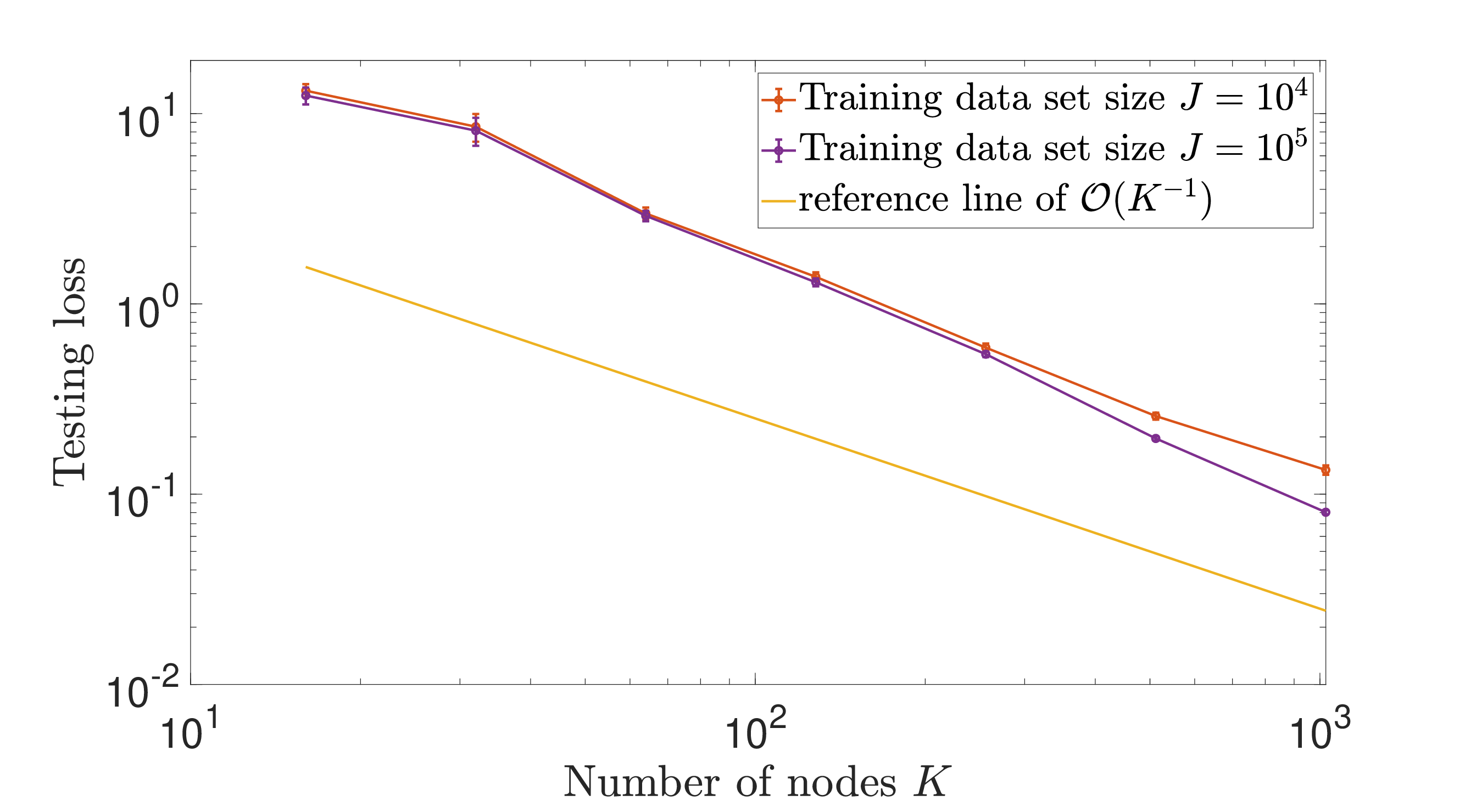}
    \caption{Gradient descent for $\omega$ }
  \end{subfigure}
  \caption{ The empirical testing loss with increasing number of nodes $K$ of the Fourier neural network, using training data set size $J=10^4$ and $J=10^5$, respectively. }
  \label{fig:general_loss_varying_J}
\end{figure}

In Figure~\ref{fig:general_loss_varying_J} we plot the empirical generalization error $\mathcal{L}_{\mathrm{emp}}$ with increasing number of nodes $K$, by implementing the adaptive Metropolis method and the gradient descent method for updating the frequency parameters $\{\omega_k\}_{k=1}^K$, respectively. To identify the effect of the training data set size $J$, for each method we also vary the value of $J$ from $10^4$ to $10^5$, and make a comparison between the corresponding empirical generalization losses. For computing $\mathcal{L}_{\mathrm{emp}}$ we choose parameters $\alpha_1=\alpha_2=1$, and employ $Q=32$ independent replicas to obtain 95\% confidence intervals of the empirical testing loss values.

For both methods, the empirical generalization error $\mathcal{L}_{\mathrm{emp}}$ are observed to decrease for increasing values of $K$, following the $\mathcal{O}(K^{-1})$ estimate in \eqref{gen_J}. Moreover, increasing the training data set size $J$ from $10^4$ to $10^5$ leads to reduced $\mathcal{L}_{\mathrm{emp}}$ loss. This declining effect in $\mathcal{L}_{\mathrm{emp}}$ becomes more pronounced for larger values of $K$, since as $K$ increases the error term $\mathcal{O}(J^{-\frac{1}{2}})$ gradually becomes comparable to $\mathcal{O}(K^{-1})$, whereas for smaller values of $K$ the $\mathcal{O}(K^{-1})$ term is still the dominating error contributor.

\begin{figure}[ht!] 
  \centering
  \begin{subfigure}{0.33\textwidth}
    \includegraphics[width=\linewidth]{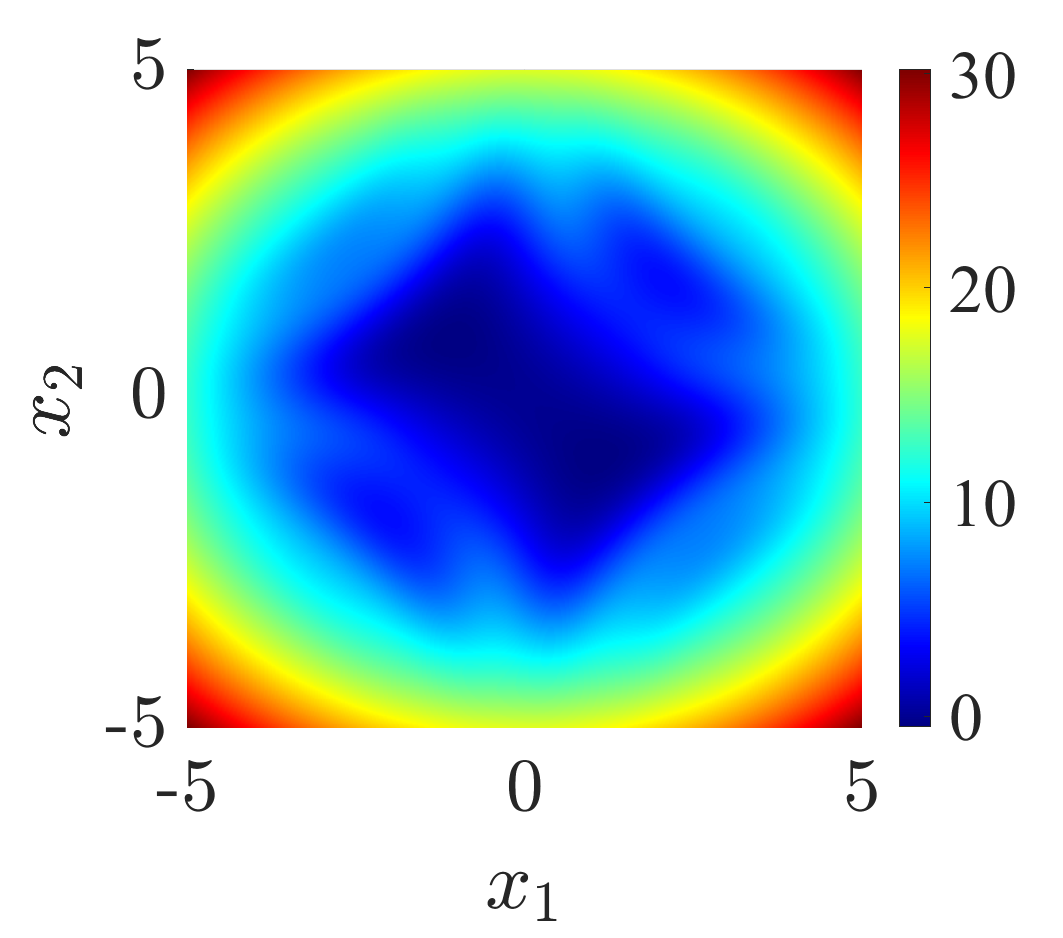}
    \caption{ $V(x)$ }
  \end{subfigure}
  \hspace{-0.02\textwidth} %
  \begin{subfigure}{0.33\textwidth}
    \includegraphics[width=\linewidth]{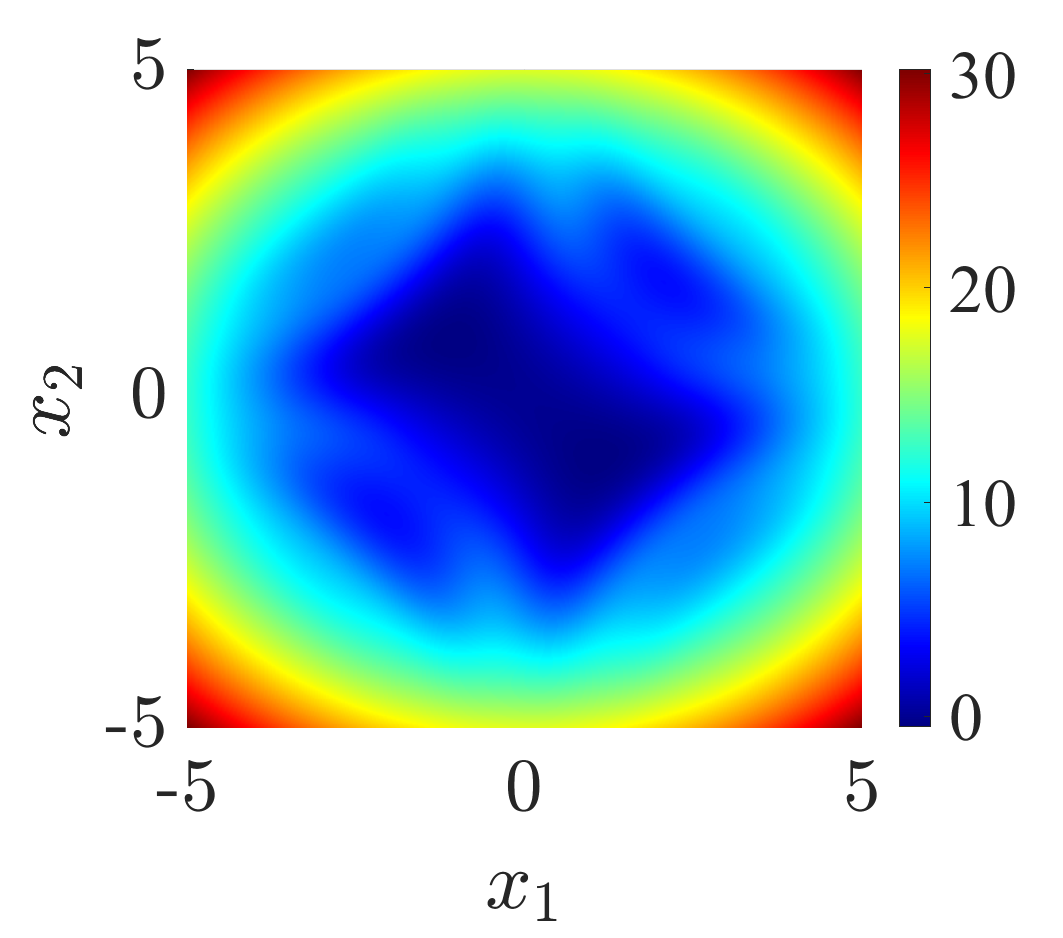}
    \caption{ $\Bar{v}_r(x)$ }
  \end{subfigure}
  \hspace{-0.02\textwidth} %
  \begin{subfigure}{0.33\textwidth}
    \includegraphics[width=\linewidth]{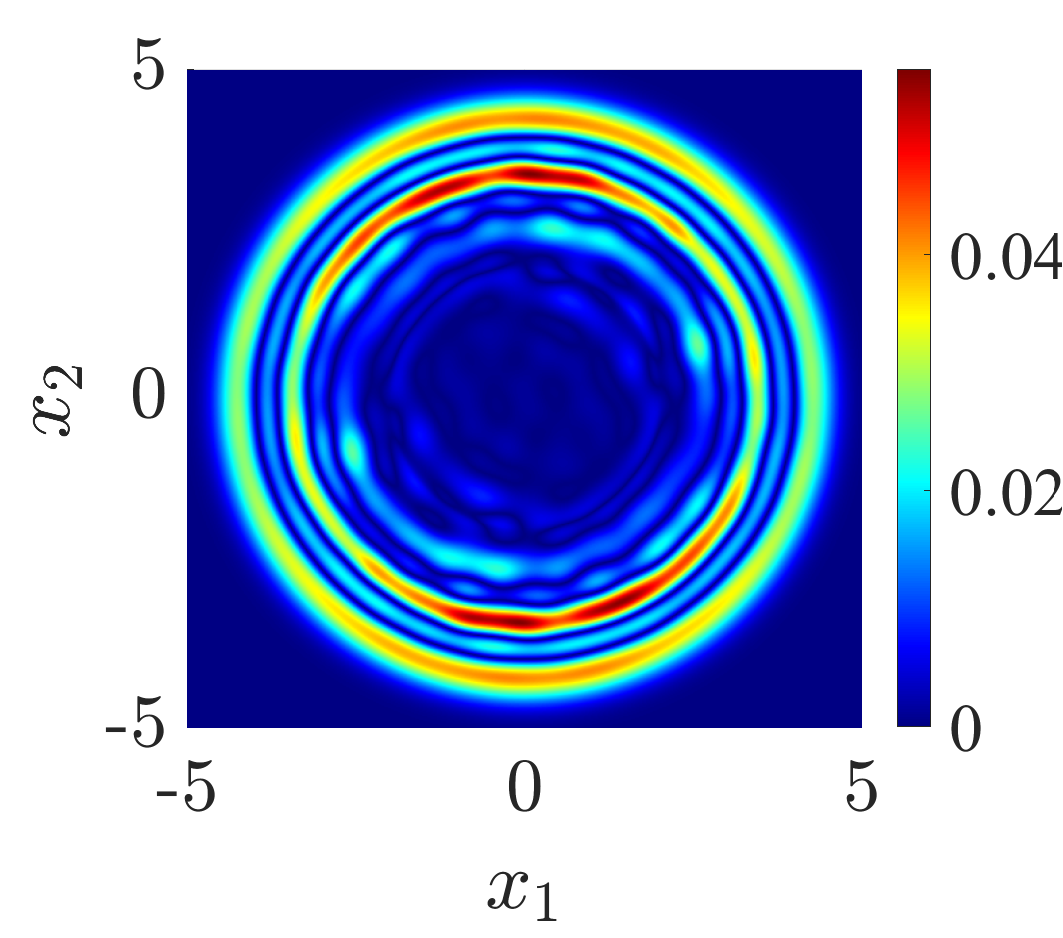}
    \caption{Pointwise difference }
  \end{subfigure}
  \caption{ Visualization of the target potential function $V(x)$, the reconstructed potential function $\Bar{v}_r(x)$ by the Fourier neural network with number of nodes $K=1024$ and data set size $J=10^5$, along with their pointwise difference. }
  \label{fig:visualize_v}
\end{figure}

To present a more direct examination of the training result, we plot the target potential function $V(x)$, the reconstructed potential function $\Bar{v}_r(x)$ using the Fourier neural network approximation $\bar{v}(x)$, along with their pointwise difference in Figure~\ref{fig:visualize_v}. Specifically $\Bar{v}_r(x)$ is obtained by
\[\Bar{v}_r(x)=\mathrm{Re}\big(\bar{v}(x)\chi_b(x)\big) +V(x)(1-\chi_b(x)),\,\]
where $\chi_b(x)=\mathbbm{1}_{\{|x|\leq R_b\}} + \mathbbm{1}_{\{|x|> R_b\}}\exp{\big( -\frac{(|x|-R_b)^2}{2\delta_b} \big)}$
with $R_b=3$ and $\delta_b=0.2$, denoting the smooth cutoff function at radial distance $R_b$. The Fourier neural network function $\bar{v}(x)$ employs $K=1024$ nodes and is trained with a data set of size $J=10^5$ applying the adaptive Metropolis algorithm.

The training data set is generated using a hybrid approach involving Langevin dynamics sampling under two inverse temperatures, namely $\beta_1=1$ and $\beta_2=0.3$, with each contributing equally to the data set. This sampling method provides a sufficient amount of data within the domain $|x|<R_b$, and facilitates a more accurate approximation therein. We further compare this mixed-temperature sampling technique with the sampling solely under a single inverse temperature, as discussed in Section~\ref{subsec_nume_one_beta}.

\subsection{The approximation of correlation observables}\label{subsec_nume_corr_func}
\if\JOURNAL1
\leavevmode \\
\fi
Based on the trained potential function $\Bar{v}(x)$ with the Fourier network representation and the corresponding reconstructed potential function $\Bar{v}_r(x)$, we calculate the approximated correlation function $\bar{\mathcal{C}}_{AB}(\tau)$ for different correlation time $\tau$ following \eqref{C_ab_NN} as
\begin{equation}\label{Corr_func_tau}
    \begin{aligned}
    \bar{\mathcal{C}}_{AB}(\tau)&=\int_{\mathbb{R}^4} A\big(\bar z_\tau(z_0)\big)B(z_0)\, \Bar{\mu}(z_0)\, \mathrm{d}z_0 \\
    &=\int_{\rtset^{2}} \int_{\rtset^{2}}A\big(\bar z_\tau(z_0)\big)B(z_0) \,\Bar{\mu}_x(x)\,\phi(p)\,\mathrm{d}x_0 \,\mathrm{d}p_0\,,
    \end{aligned}
\end{equation}
where
\[
\bar \mu(z)= \frac{e^{-\beta(|p|^2/2+\bar v_r(x))}}{\int_{\rtset^{4}}  e^{-\beta(|p|^2/2+\bar v_r(x))}\,\mathrm{d}z},\quad \Bar{\mu}_x(x)=\frac{e^{-\beta \Bar{v}_r(x)}}{\int_{\rset^{2}}e^{-\beta \Bar{v}_r(x)}\,\mathrm{d}x}\,,
\]
and $\phi(p)$ is the probability density function of multivariate normal distribution $\mathcal{N}(0, \frac{1}{\beta}\mathrm{I}_{2\times 2})$. The formulation \eqref{Corr_func_tau} admits a convenient Monte Carlo approximation based on independent samples of the initial states $z_0^{(m)}=(x_0^{(m)},p_0^{(m)})$ for $m=1,2,\dots, M$ in the phase space as
\begin{equation}\label{Corr_func_Monte_Carlo}
    \bar{\mathcal{C}}_{AB}(\tau)\simeq \frac{1}{M}\sum_{m=1}^M A(\Bar{z}_\tau(z_0^{(m)})\,B(z_0^{(m)})\,,
\end{equation}
where $\{x_0^{(m)} \}_{m=1}^M$ is the set of initial position samples drawn by the overdamped Langevin dynamics \eqref{langevin} using the reconstructed potential $\Bar{v}_r$, and $\{p_0^{(m)} \}_{m=1}^M$ is the set of initial momentum samples following $\mathcal{N}(0, \frac{1}{\beta}\mathrm{I}_{2\times 2})$ distribution. The approximated auto-correlation function $\bar{\mathcal{C}}_{AB}(\tau)$
is then compared with its reference value $\mathcal{C}_{AB}(\tau)$ 
which is evaluated by a Monte Carlo approximation similar to \eqref{Corr_func_Monte_Carlo} as
\begin{equation}\label{Corr_ref_Monte_Carlo}
\begin{aligned}
\mathcal{C}_{AB}(\tau) &=\int_{\rtset^{2}} \int_{\rtset^{2}}A\big( z_\tau(z_0)\big)B(z_0) \,\mu_x(x)\,\phi(p)\,\mathrm{d}x_0 \,\mathrm{d}p_0\\
&\simeq \frac{1}{M_{\mathrm{ref}}}\sum_{m=1}^{M_{\mathrm{ref}} }A(z_\tau(z_0^{(m)})\,B(z_0^{(m)})\,,
\end{aligned}
\end{equation}
 where $\mu_x(x)=e^{-\beta V(x)}/{\int_{\rset^{2}}e^{-\beta V(x)}\,\mathrm{d}x}$ is the Gibbs density associated with 
 the target potential $V(x)$, while the initial momentum samples are drawn from the same multivariate normal distribution $\mathcal{N}(0, \frac{1}{\beta}\mathrm{I}_{2\times 2})$ as in \eqref{Corr_func_Monte_Carlo}.

For the numerical validation we compute the approximated auto-correlation function $\bar{\mathcal{C}}_{x_1,x_1}(\tau)$ between the first component of the position observables at time $t=0$ and $t=\tau$, by taking $A\big(\bar z_\tau(z_0)\big)=\Bar{x}_1(\tau)$ and $B(z_0)=x_1(0)$ in \eqref{Corr_func_Monte_Carlo}, where
\[
\bar z_t=[\,\Bar{x}_t,\Bar{p}_t\,]^T=[ \,\Bar{x}_1(t),\,\Bar{x}_2(t),\,\Bar{p}_1(t),\,\Bar{p}_2(t)\,]^T\,,
\]
satisfying the dynamics
\begin{equation}\label{dz_bar_dt}
\dot{\bar z}_t = \bar f(\bar z_t)=[\,\Bar{p}_1(t),\, \Bar{p}_2(t),\,-\partial_{x_1}\,\Bar{v}_r(\Bar{x}_t), \,-\partial_{x_2}\,\Bar{v}_r(\Bar{x}_t)\,]^T\,.
\end{equation}
Here we use the notation $T$ to represent the transpose of a row vector. Similarly we also calculate the approximation of the auto-correlation function $\bar{\mathcal{C}}_{p_1,p_1}(\tau)$ between the first components $\Bar{p}_1(\tau)$ and $p_1(0)$ of the momentum observables using \eqref{Corr_func_Monte_Carlo}. 

\begin{figure}%
\centering
\includegraphics[height=0.5\textwidth]{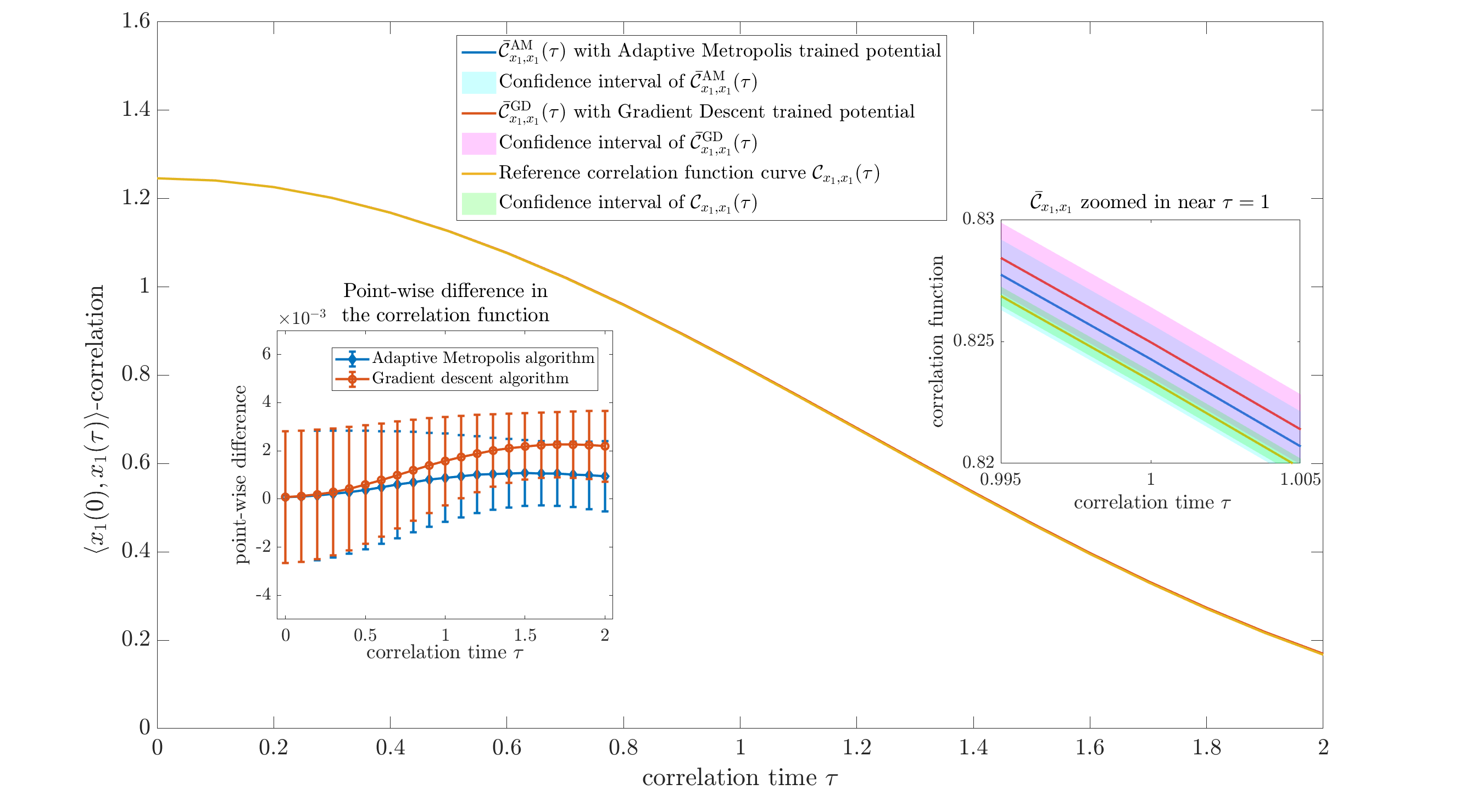}
\caption{ Approximation of the position auto-correlation function $\Bar{\mathcal{C}}_{x_1,x_1}(\tau)$, using Fourier neural network with $K=1024$ nodes, data set size $J=10^{5}$, $\beta=1$, and Monte Carlo sample size $M=2^{21}$. }
\label{fig:Corr_func_x}
\end{figure}

The reference values $\mathcal{C}_{x_1,x_1}(\tau)$ and $\mathcal{C}_{p_1,p_1}(\tau)$ corresponding to $\bar {\mathcal{C}}_{x_1,x_1}(\tau)$ and $\bar {\mathcal{C}}_{p_1,p_1}(\tau)$ are obtained by taking $A\big(z_\tau(z_0)\big)=x_1(\tau)$,  $B(z_0)=x_1(0)$ and $A\big(z_\tau(z_0)\big)=p_1(\tau)$,  $B(z_0)=p_1(0)$ in \eqref{Corr_ref_Monte_Carlo} respectively, where
\[
 z_t=[\,x_t,p_t\,]^T=[ \,x_1(t),\,x_2(t),\,p_1(t),\,p_2(t)\,]^T\,,
\]
has the dynamics described by 
\begin{equation}\label{dz_dt}
\dot{z}_t = f( z_t)=[\,p_1(t),\, p_2(t),\,-\partial_{x_1}\,V(x_t), \,-\partial_{x_2}\,V(x_t)\,]^T\,.
\end{equation}
In Figure ~\ref{fig:Corr_func_x}, we present the approximated auto-correlation function curves $\bar{\mathcal{C}}_{x_1,x_1}^{\mathrm{AM}}(\tau)$ and  $\bar{\mathcal{C}}_{x_1,x_1}^{\mathrm{GD}}(\tau)$ at the inverse temperature $\beta=1$, generated by Fourier neural networks trained with adaptive Metropolis method and gradient descent method, respectively. These curves are juxtaposed with the reference curve $\mathcal{C}_{x_1,x_1}(\tau)$. Given the close proximity of these curves, we magnify the plot around $\tau=1$ and include an inset depicting the pointwise differences $\bar{\mathcal{C}}_{x_1,x_1}^{\mathrm{AM}}(\tau_i)-\mathcal{C}_{x_1,x_1}(\tau_i)$ and $\bar{\mathcal{C}}_{x_1,x_1}^{\mathrm{GD}}(\tau_i)-\mathcal{C}_{x_1,x_1}(\tau_i)$ for $i=0,1,\dots,20$, up to the correlation time $\tau=2$, where $\tau_i = i \Delta\tau$ with discretization in correlation time $\Delta \tau = 0.1 $. 

\begin{figure}%
\centering
\includegraphics[height=0.5\textwidth]{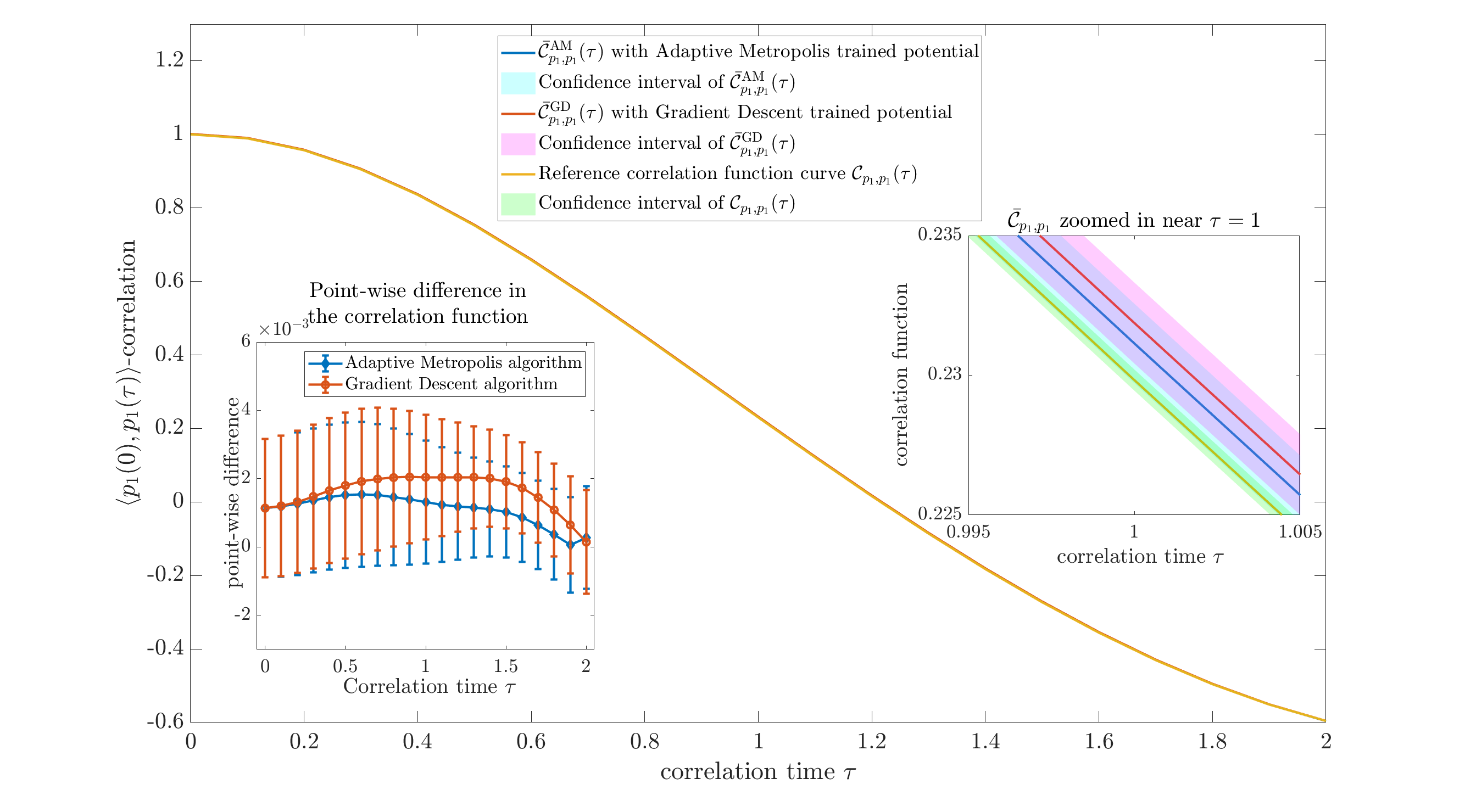}
\caption{ Approximation of the momentum auto-correlation function $\Bar{\mathcal{C}}_{p_1,p_1}(\tau)$, using Fourier neural network with $K=1024$ nodes, data set size $J=10^{5}$, $\beta=1$, and Monte Carlo sample size $M=2^{21}$. }
\label{fig:Corr_func_p}
\end{figure}

For the Monte Carlo approximation of each point in the auto-correlation function $\bar{\mathcal{C}}_{x_1,x_1}(\tau_i)$, we compute the corresponding 95\% confidence interval based on $Q=32$ independent replicas, employing a total sample size $M = 2^{21}$. For the reference values $\mathcal{C}_{x_1,x_1}(\tau_i)$, we still use $Q=32$ independent replicas to evaluate the confidence intervals, while the total Monte Carlo sample size is increased to $M_{\mathrm{ref}}=2^{26}$ in order to acquire a reference curve with high accuracy. The dynamics \eqref{dz_bar_dt} and \eqref{dz_dt} are numerically approximated by implementing the 2nd-order velocity Verlet method \cite{Leimkuhler_Matthews_MD_book} with time step size $\Delta t=5\times 10^{-3}$, where the associated time discretization error is numerically tested to be one order of magnitude smaller than the statistical error by the Monte Carlo estimation. A similar plot for the momentum auto-correlation function $\bar{\mathcal{C}}_{p_1,p_1}(\tau)$ is given in Figure~\ref{fig:Corr_func_p}, using the same parameter settings as in Figure~\ref{fig:Corr_func_x}.

From Figure~\ref{fig:Corr_func_x} we observe that the pointwise difference gradually rises as the correlation time $\tau$ increases, while the confidence intervals of all three curves overlap in both figures~\ref{fig:Corr_func_x} and ~\ref{fig:Corr_func_p}. To further survey the performance of the  reconstructed potential $\Bar{v}_r(x)$, 
we plot the $L^1$-difference between the curves $\|\mathcal{C}_{x_1,x_1}(\tau)-\bar{\mathcal{C}}_{x_1,x_1}(\tau)\|_{L^1}$ and $\|\mathcal{C}_{p_1,p_1}(\tau)-\bar{\mathcal{C}}_{p_1,p_1}(\tau)\|_{L^1}$, which are approximately evaluated by
\[
\begin{aligned}
\|\mathcal{C}_{x_1,x_1}(\tau)-\bar{\mathcal{C}}_{x_1,x_1}(\tau)\|_{L^1}&\approx \sum_{i=0}^{20} \big|\mathcal{C}_{x_1,x_1}(\tau_i)-\bar{\mathcal{C}}_{x_1,x_1}(\tau_i)\big|\Delta \tau\,,\\
\mbox{ and}\quad \|\mathcal{C}_{p_1,p_1}(\tau)-\bar{\mathcal{C}}_{p_1,p_1}(\tau)\|_{L^1}&\approx \sum_{i=0}^{20} \big|\mathcal{C}_{p_1,p_1}(\tau_i)-\bar{\mathcal{C}}_{p_1,p_1}(\tau_i)\big| \Delta \tau\,.
\end{aligned}
\]

\begin{figure}[ht!] 
  \centering
  \begin{subfigure}{0.51\textwidth}
    \includegraphics[width=\linewidth]{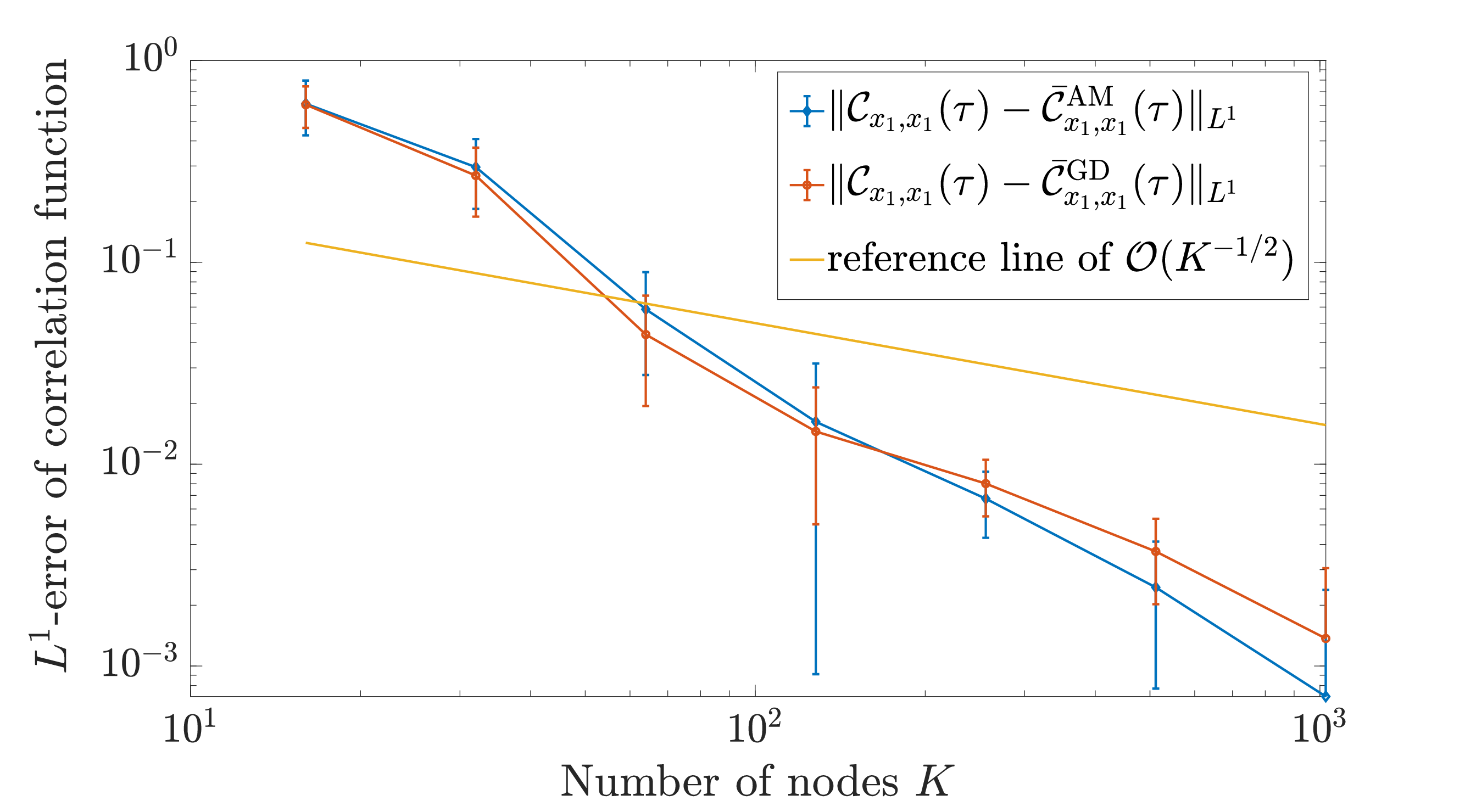}
    \caption{ $L^1$-error of $\bar{\mathcal{C}}_{x_1,x_1}(\tau)$  }
  \end{subfigure}
  \hspace{-0.04\textwidth} %
  \begin{subfigure}{0.51\textwidth}
    \includegraphics[width=\linewidth]{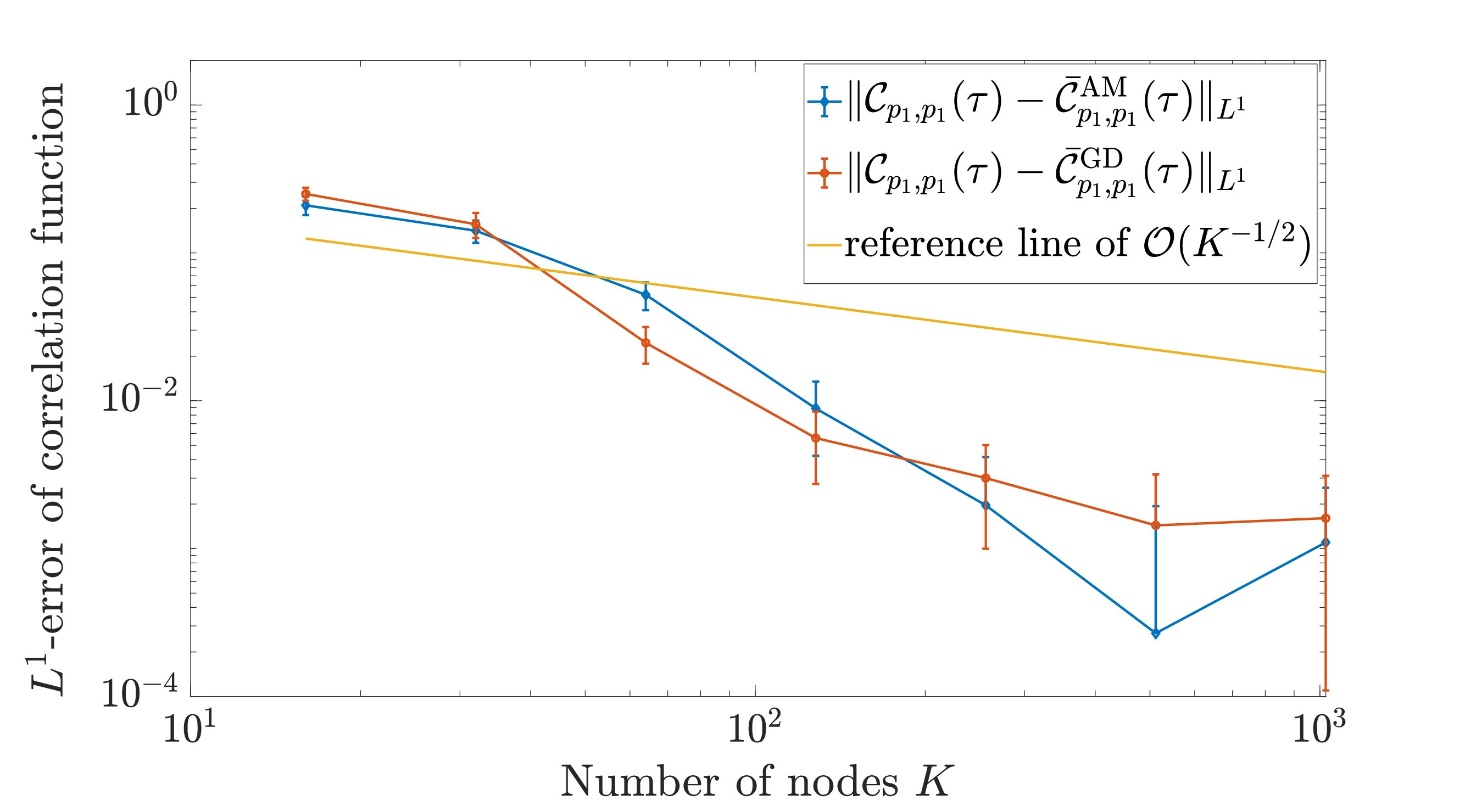}
    \caption{$L^1$-error of $\bar{\mathcal{C}}_{p_1,p_1}(\tau)$ }
  \end{subfigure}
  \caption{ The $L^1$-difference of the approximated auto-correlation function $\|\mathcal{C}_{x_1,x_1}(\tau)-\bar{\mathcal{C}}_{x_1,x_1}(\tau)\|_{L^1}$ and $\|\mathcal{C}_{p_1,p_1}(\tau)-\bar{\mathcal{C}}_{p_1,p_1}(\tau)\|_{L^1}$, with increasing number of nodes $K$ in the Fourier neural network. The training data set size $J=10^5$, and the statistical uncertainties are evaluated with $Q=32$ independent replicas.}
  \label{fig:Corr_L1_diff}
\end{figure}

Corroborating the Theorem~\ref{thm}, Figure~\ref{fig:Corr_L1_diff} numerically confirms a decreasing trend in the $L^1$-error of the approximated correlation functions $\|\mathcal{C}_{x_1,x_1}(\tau)-\bar{\mathcal{C}}_{x_1,x_1}(\tau)\|_{L^1}$ and $\|\mathcal{C}_{p_1,p_1}(\tau)-\bar{\mathcal{C}}_{p_1,p_1}(\tau)\|_{L^1}$ with an increase in the number of nodes $K$. As $K$ grows, the declining trend of the the $L^1$-error gradually aligns with the analytical error bound $\mathcal{O}(K^{-1/2})$ provided in \eqref{C_b}. Since the statistical uncertainty arising from the Monte Carlo estimation and the influence of the data set size $J$ emerges as the predominant error, it remains challenging to precisely evaluate the convergence rate with respect to $K$.

\subsection{Training data using different temperatures}\label{subsec_nume_one_beta}
\if\JOURNAL1
\leavevmode \\
\fi
Throughout the numerical experiments of Sections~\ref{subsec_nume_general_error} and~\ref{subsec_nume_corr_func}, we employ a hybrid sampling technique with overdamped Langevin dynamics \eqref{langevin} under two fixed inverse temperatures $\beta_1=1$ and $\beta_2=0.3$, with each dynamics contributing equally to the final training data set. 

\begin{figure}[ht!] 
  \centering
  \begin{subfigure}{0.51\textwidth}
    \includegraphics[width=\linewidth]{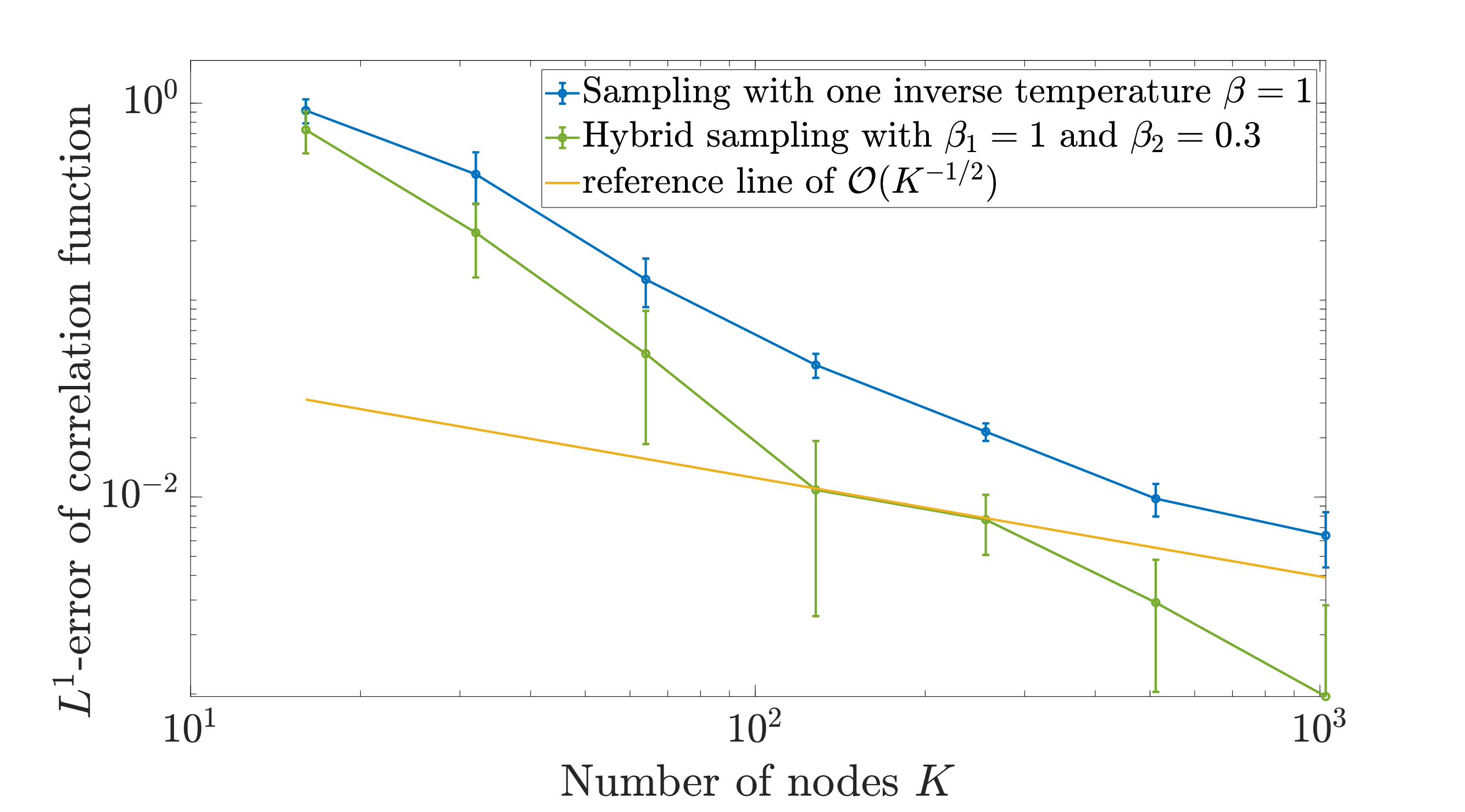}
    \caption{ $L^1$-error of $\bar{\mathcal{C}}_{x_1,x_1}^{\mathrm{GD}}(\tau)$  }
  \end{subfigure}
  \hspace{-0.04\textwidth} %
  \begin{subfigure}{0.51\textwidth}
    \includegraphics[width=\linewidth]{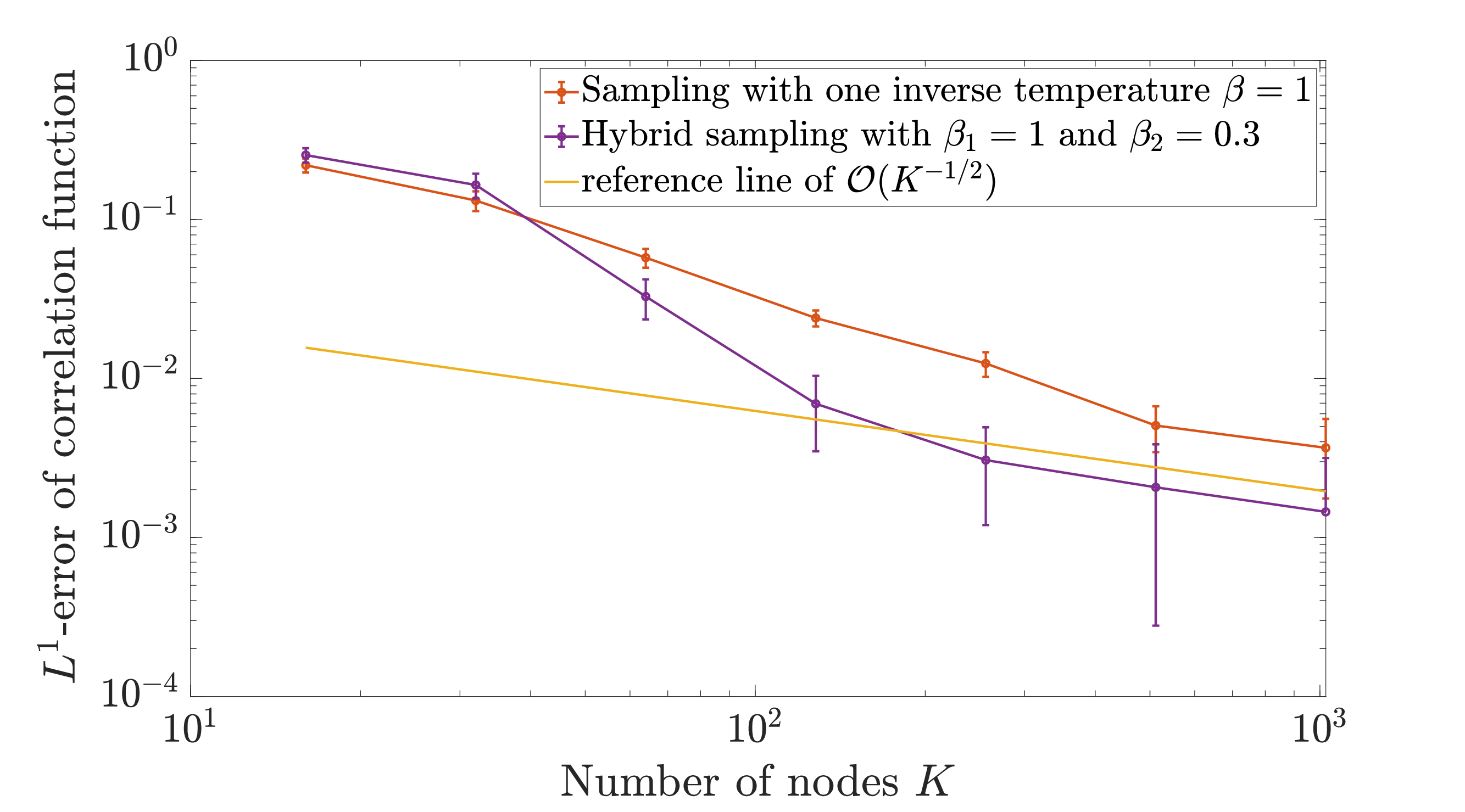}
    \caption{$L^1$-error of $\bar{\mathcal{C}}_{p_1,p_1}^{\mathrm{GD}}(\tau)$ }
  \end{subfigure}
  \caption{ The $L^1$-error of the approximated auto-correlation functions $\|\mathcal{C}_{x_1,x_1}(\tau)-\bar{\mathcal{C}}_{x_1,x_1}^{\mathrm{GD}}(\tau)\|_{L^1}$ and $\|\mathcal{C}_{p_1,p_1}(\tau)-\bar{\mathcal{C}}_{p_1,p_1}^{\mathrm{GD}}(\tau)\|_{L^1}$, with increasing number of nodes $K$ in the Fourier neural network trained with gradient descent method on the frequencies. The training data sets are of size $J=10^4$ and are sampled with two different strategies. The statistical uncertainties are evaluated with $Q=32$ independent replicas.}
  \label{fig:compare_sampling_corr}
\end{figure}

We choose  $\beta_1=1$ to guarantee sufficient amount of relevant data sampled at the same temperature as that for evaluating the approximated correlation observables. Additionally, selecting $\beta_2 = 0.3$ facilitates sampling at a higher temperature, increasing the probability of sampled data points overcoming potential barriers. Consequently, more data points explore a larger phase space domain, reducing the confinement near the minimum of the potential surface. The combination of two sampling dynamics yields a more comprehensive training dataset. In contrast, we implement overdamped Langevin dynamics sampling under a single fixed inverse temperature $\beta = 1$, and compare its results with the hybrid sampling strategy in Figure~\ref{fig:compare_sampling_corr}.

Upon comparing the $L^1$-error of the auto-correlation function approximations resulting from the two sampling strategies in Figure~\ref{fig:compare_sampling_corr}, both exhibit an asymptotic decaying trend consistent with the $\mathcal{O}(K^{-1/2})$ estimate in \eqref{C_b}. Notably, for large values of $K$, the $L^1$-error incurred by the hybrid sampling strategy (depicted by the green curve in the left panel and the purple curve in the right panel) is significantly smaller than that yielded by sampling with a single inverse temperature (represented by the blue curve in the left panel and the red curve in the right panel), indicating a clear improvement offered by the hybrid sampling approach.

\subsection{%
Optimization of 
the regularized loss function 
$\mathcal{L}_{\mathrm{R}}$ 
}\label{subsec:optimization_alg}
\if\JOURNAL1
\leavevmode \\
\fi
To optimize the amplitude coefficients $\eta$ 
for a fixed set of frequency parameters $\omega$, we explore three variants of the objective function $\mathcal{L}_{\mathrm{R}}$ in \eqref{regularized_loss}, namely,
\[
\begin{aligned}
    \mathcal{L}_{\mathrm{R}_1}(\eta) &:=
\frac{1}{J}\sum_{j=1}^J\big(
\alpha_1|v(x_j)-\bar{v}(x_j;\eta)|^2 + \alpha_2|v'(x_j)-\bar{v}'(x_j;\eta)|^2\big) 
+\lambda_1\sum_{k=1}^K|\eta_k|^2\,,\\
  \mathcal{L}_{\mathrm{R}_2}(\eta) &:=
\mathcal{L}_{\mathrm{R}_1}(\eta)+\lambda_2(\sum_{k=1}^K|\eta_k|^2)^2\,,\\ %
 \mathcal{L}_{\mathrm{R}_3}(\eta) &:=
\mathcal{L}_{\mathrm{R}_2}(\eta)+\lambda_3\max(\sum_{k=1}^K|\eta_k|-\tilde{C},0) \,,
\end{aligned}
\]
where $\mathcal{L}_{\mathrm{R}_1}(\eta)$ is the standard squared training loss computed from the discrepancy between the network output function $\bar{v}$ and the training data $v$, with Tikhonov regularization $\lambda_1\sum_{k=1}^K|\eta_k|^2$. In $\mathcal{L}_{\mathrm{R}_2}(\eta)$ and $\mathcal{L}_{\mathrm{R}_3}(\eta)$, we augment the objective function with a regularization term 
\if\JOURNAL1
\\
\fi
$\lambda_2(\sum_{k=1}^K|\eta_k|^2)^2$. Additionally, for $\mathcal{L}_{\mathrm{R}_3}(\eta)$, an extra penalization term $\lambda_3\max\big(\sum_{k=1}^K|\eta_k|-\tilde{C},0\big)$ is introduced. 

\begin{figure}[ht!] 
  \centering
  \begin{subfigure}{0.51\textwidth}
    \includegraphics[width=\linewidth]{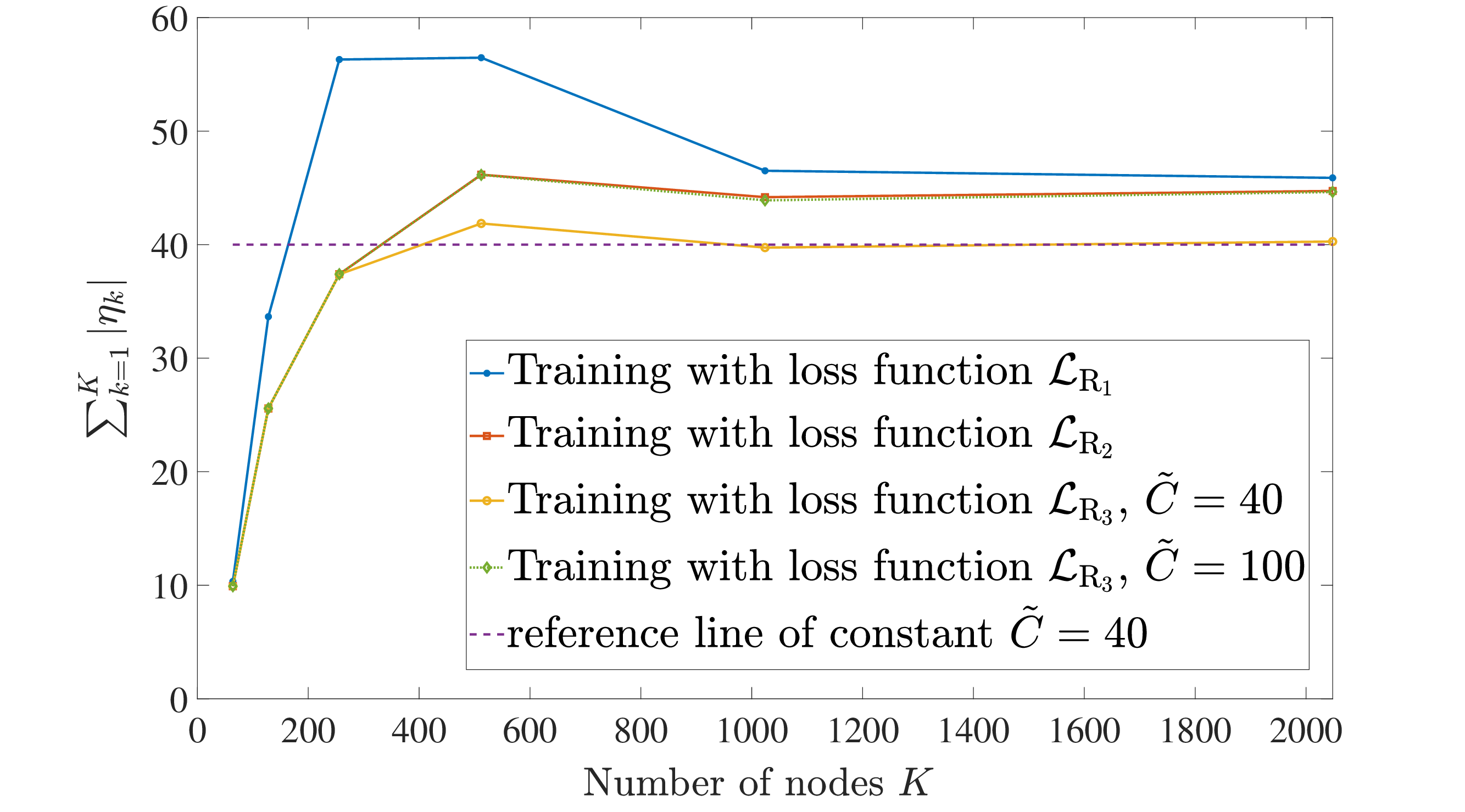}
    \caption{  $\sum_{k=1}^K|\eta_k|$ }
    \label{fig:training_eta_sum_subfig1}
  \end{subfigure}
  \hspace{-0.04\textwidth} %
  \begin{subfigure}{0.51\textwidth}
    \includegraphics[width=\linewidth]{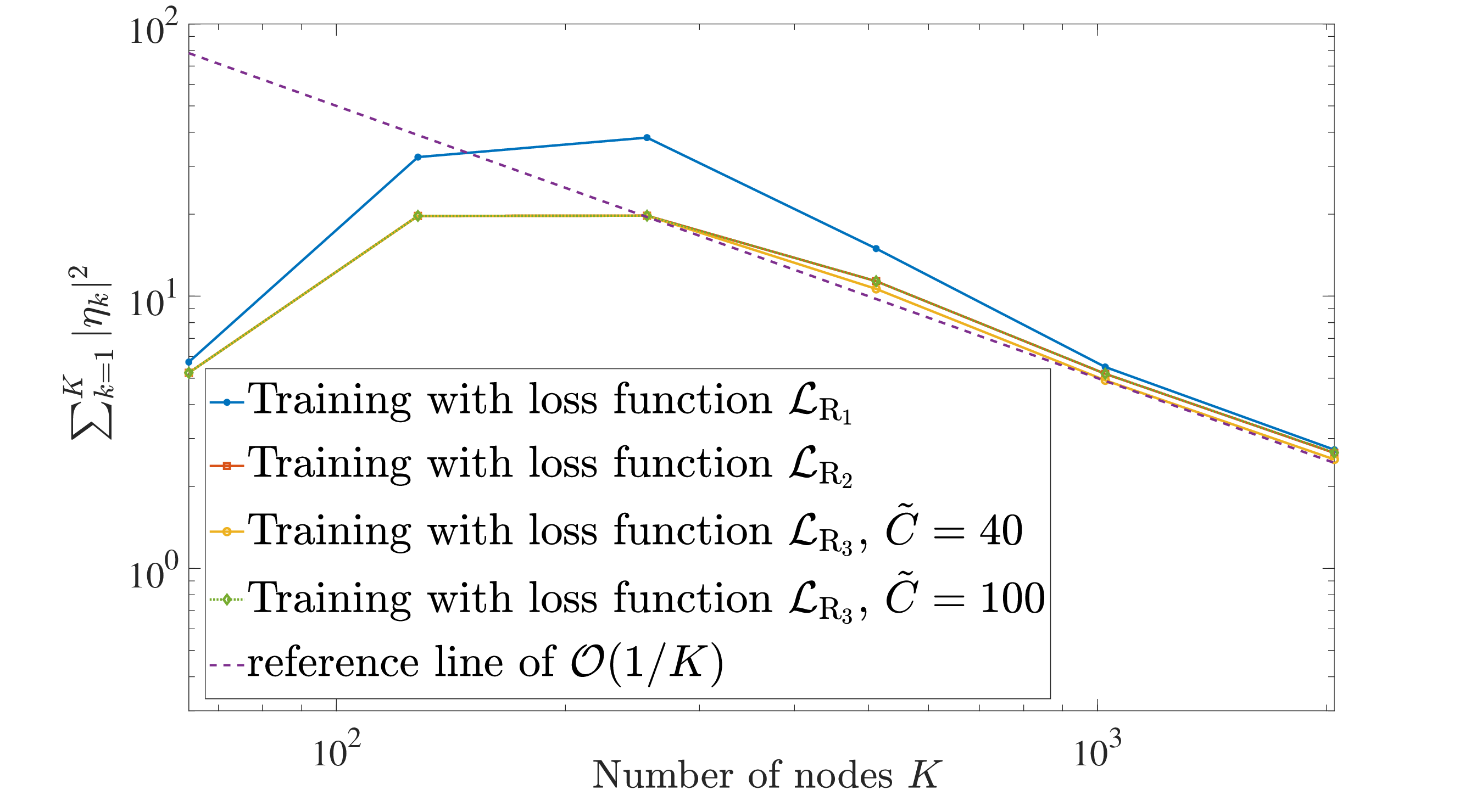}
    \caption{$\sum_{k=1}^K|\eta_k|^2$}
    \label{fig:training_eta_sum_subfig2}
  \end{subfigure}
  \caption{ The sum of absolute amplitudes and squared amplitudes, obtained by training with three variants of regularized loss function $\mathcal{L}_{\mathrm{R}_1}(\eta)$, $\mathcal{L}_{\mathrm{R}_2}(\eta)$, and $\mathcal{L}_{\mathrm{R}_3}(\eta)$, using parameters $\alpha_1=1$, $\alpha_2=1$, $\lambda_1 = 10^{-2}$, $\lambda_2=10^{-3}$, $\lambda_3 = 10^{-2}$, and $\tilde{C}=40$ or $100$.}
  \label{fig:training_eta_sum}
\end{figure}

In response to the optimal amplitude coefficients in \eqref{eta_F} with frequency parameters sampled from the optimal density $\rho_\ast$ given by \eqref{rho_star_optimal}, we have
\begin{equation}\label{optimal_eta_sum}
\hat\eta_k =\frac{1}{K} \,\frac{\widehat{v}(\omega_k)}{\PN_\ast(\omega_k)}\,,\ k=1,\ldots, K\,,\ \Rightarrow\ \sum_{k=1}^K|\Hat{\eta}_k|= \mathcal{O}(1)\,, \ \ \sum_{k=1}^K|\Hat{\eta}_k|^2=\mathcal{O}(K^{-1})\,.
\end{equation}
We evaluate the empirical sum of the absolute amplitudes $\sum_{k=1}^K|\eta_k|$ and the sum of squared amplitudes $\sum_{k=1}^K|\eta_k|^2$ yielded by the above three variants of regularized loss functions, and depict the results in Figure~\ref{fig:training_eta_sum}

It is well known that imposing a size constraint on the amplitude coefficients through Tikhonov regularization effectively alleviates the problem of cancellation between positive and negative coefficients with large absolute values \cite{elements_of_stat_learning}.  This property is demonstrated in Figure~\ref{fig:training_eta_sum}, where all the curves in Subfigure~\ref{fig:training_eta_sum_subfig1} indicate a bounded sum of absolute amplitudes, and the sum of squared amplitudes in Subfigure~\ref{fig:training_eta_sum_subfig2} are asymptotically consistent with the $\mathcal{O}(K^{-1})$ estimate given by \eqref{optimal_eta_sum}.

In Figure~\ref{fig:training_eta_sum}, by comparing the solid blue curve which employs the loss function $\mathcal{L}_{\mathrm{R}_1}$ with the other three curves, we observe that the introduction of the penalty term  $\lambda_2(\sum_{k=1}^K|\eta_k|^2)^2$ leads to a more effective confinement on the sum of absolute and squared amplitudes. Furthermore, by contrasting the yellow solid curve with the green dotted curve in Figure~\ref{fig:training_eta_sum_subfig1}, we note the impact of selecting a smaller parameter $\tilde{C}=40$ compared to $\tilde{C}=100$, which brings the sum of absolute amplitudes closer to $40$ as a result of the penalty term $\lambda_3\max\big(\sum_{k=1}^K|\eta_k|-\tilde{C},0\big)$.

\begin{figure}[ht!] 
  \centering
  \begin{subfigure}{0.4\textwidth}
    \includegraphics[width=\linewidth]{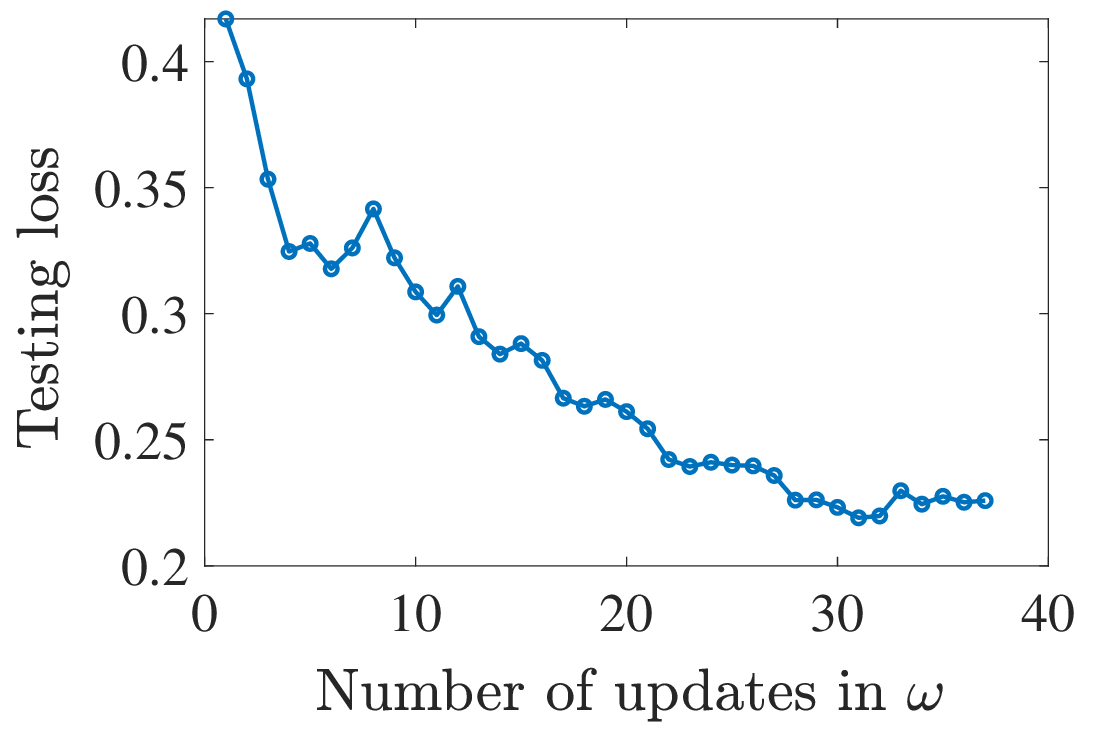}
    \caption{ adaptive Metropolis method }
    \label{fig:test_loss_omega_subfig1}
  \end{subfigure}
  \hspace{0.02\textwidth} %
  \begin{subfigure}{0.4\textwidth}
    \includegraphics[width=\linewidth]{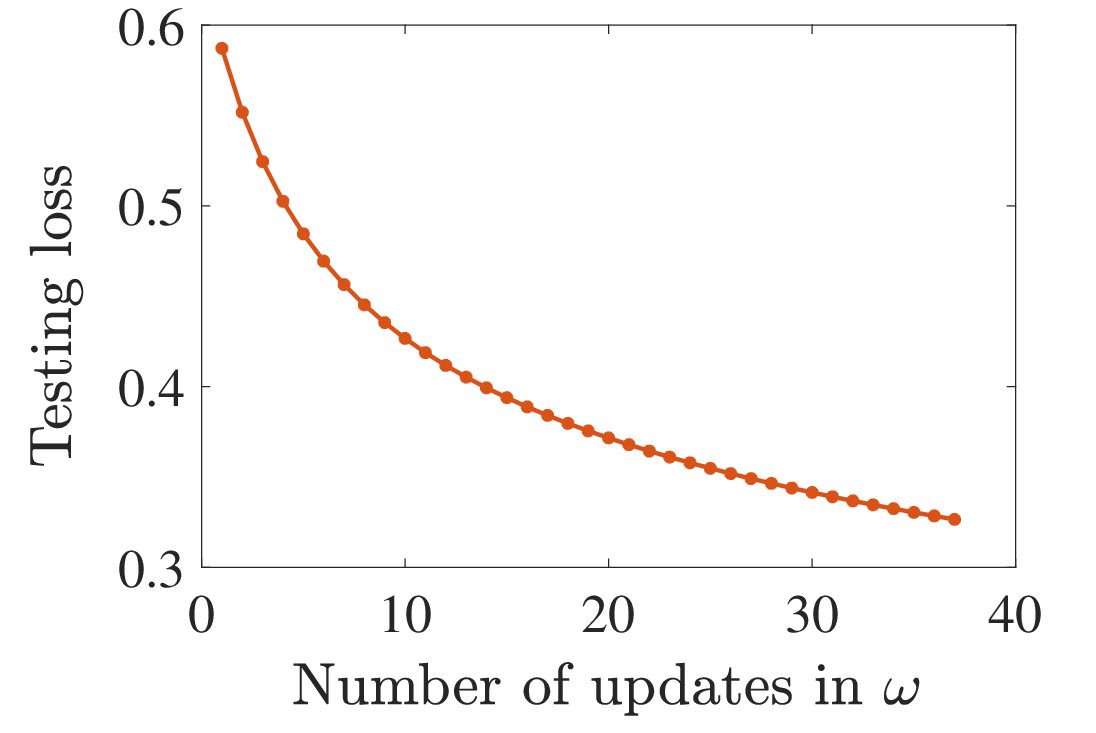}
    \caption{gradient descent method}
    \label{fig:test_loss_omega_subfig2}
  \end{subfigure}
  \caption{ The testing loss with updates in the frequency parameters $\omega$, applying adaptive Metropolis method and gradient descent method respectively.}
  \label{fig:test_loss_omega}
\end{figure}

To evaluate the effectiveness of implementing the adaptive Metropolis method and the gradient descent method for updating the frequency parameters $\omega$, we depict the testing loss against the number of iterative steps for both methods in Figure~\ref{fig:test_loss_omega}. Both methods exhibit a typical descending trend in the testing loss as the number of iterations on the frequency parameters $\omega$ increases.

\begin{figure}[ht!] 
  \centering
  \begin{subfigure}{0.33\textwidth}
    \includegraphics[width=\linewidth]{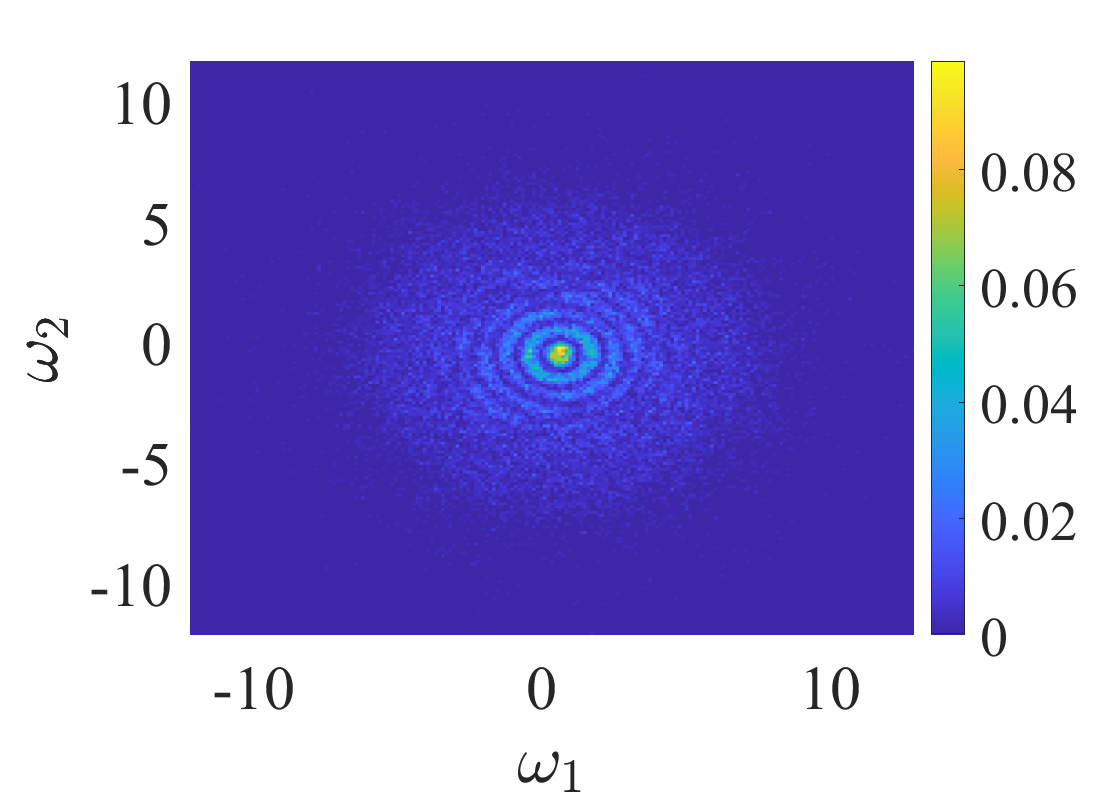}
    \caption{ Adaptive Metropolis }
  \end{subfigure}
  \hspace{-0.02\textwidth} %
  \begin{subfigure}{0.33\textwidth}
    \includegraphics[width=\linewidth]{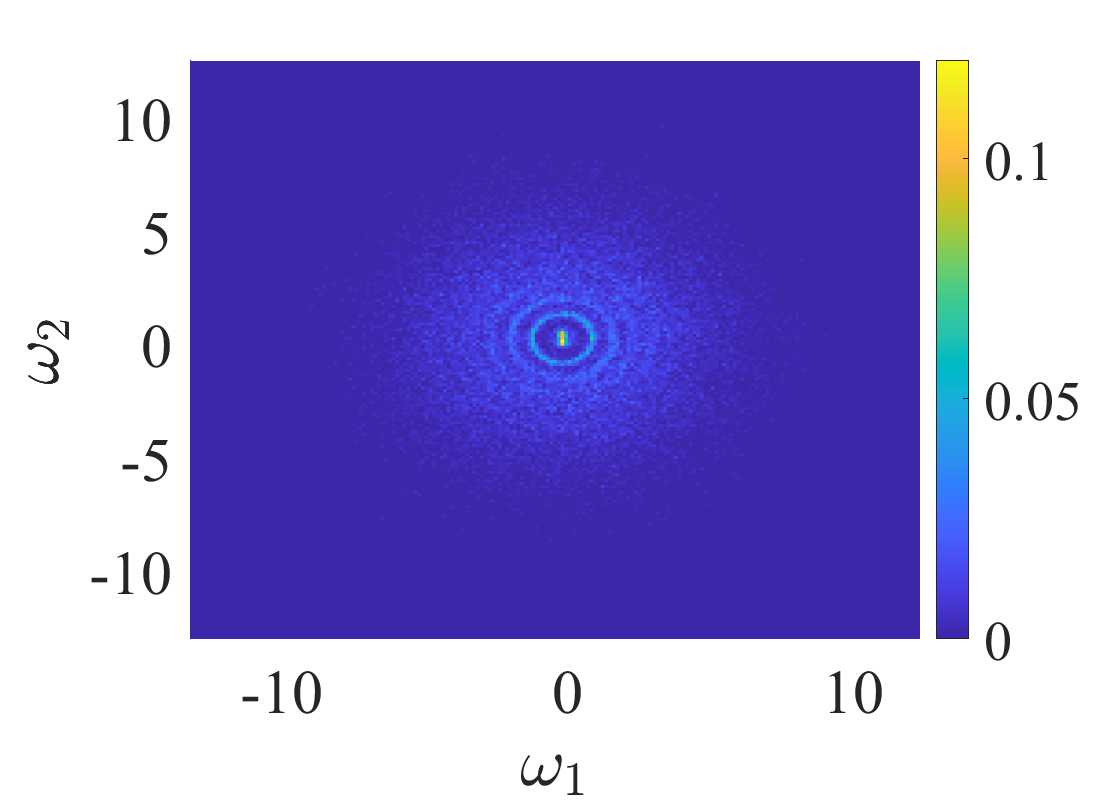}
    \caption{ Gradient Descent }
  \end{subfigure}
  \hspace{-0.02\textwidth} %
  \begin{subfigure}{0.33\textwidth}
    \includegraphics[width=\linewidth]{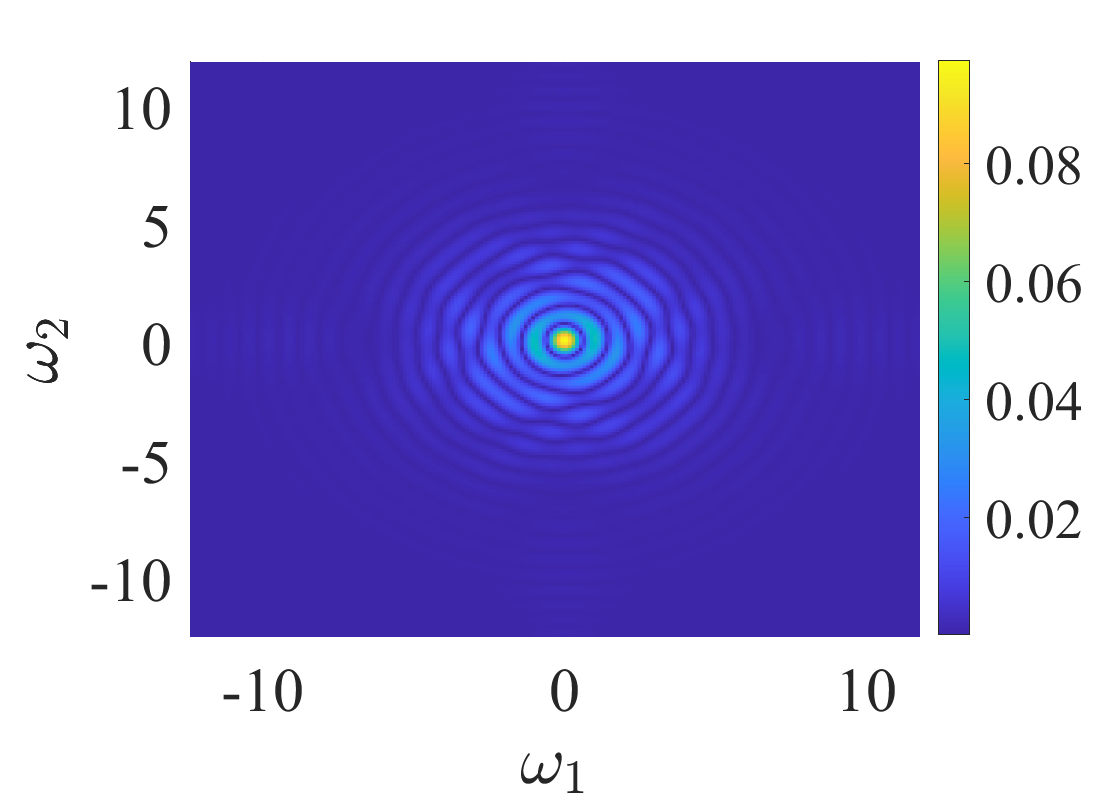}
    \caption{Optimal density $\PN_\ast(\omega)$ }
  \end{subfigure}
  \caption{ The empirical probability density function of frequencies $\omega$ obtained with the adaptive Metropolis method and the gradient descent method, and the analytically optimal density function $\PN_\ast(\omega)$. }
  \label{fig:omega_density}
\end{figure}

To further examine the quality of the trained frequency parameters $\omega$, we plot the empirical probability density function of the obtained frequency samples with $K=1024$ nodes in the Fourier neural network, by running $Q=32$ independent replicas of training. The empirical density functions by applying the adaptive Metropolis method and the gradient descent method are plotted in Figure~\ref{fig:omega_density}, along with the analytically optimal density function $\PN_\ast(\omega)$ given by \eqref{rho_star_optimal}.

Both the adaptive Metropolis and the gradient descent methods tend to sample more densely in the low-frequency domain, aligning with the theoretical optimal density function $\PN_\ast(\omega)$. However, the adaptive Metropolis method  surpasses the gradient descent method in accurately reproducing both the pattern and the scale of $\PN_\ast(\omega)$ for $\omega\in\mathbb{R}^2$.

 \section{Proof of the main theorems}\label{sec_proof} 
 We collect in this section the proofs of the main Theorem \ref{thm}, which establishes the error bound for the approximated correlation observables, and Theorem \ref{thm:generalization}, which provides an estimate of the generalization error for random Fourier feature network representation.

\subsection{Proof of Theorem \ref{thm:generalization}}\label{proof_thm_2_1}
\if\JOURNAL1
\leavevmode \\
\fi
The proof has four steps:
\begin{description}
\item[\it Step 0] Formulates the training optimization and introduces notation.
\item[\it Step 1] Estimates the training error and formulates the generalization error.
\item[\it Step 2] Derives a representation of the generalization error, 
using the regularization terms and the independence of the data points $x_j$ and the $\omega_k$ samples.
\item[\it Step 3] Combines Steps 1 and 2 to obtain \eqref{gen_J}. 
\end{description}

\begin{proof}
{\it Step 0}
We have a finite number of data points
\[\{\big(x_j,v(x_j),v'(x_j)\big)\ |\ j=1,\ldots,J\}\,,\]
where $x_j$ are independent samples from $\mu_x$, 
and  $\bar{v}(x)=\sum_{k=1}^K\eta_ke^{\mathrm{i}\omega_k\cdot x}$ is a solution to the optimization problem \eqref{opt_J}.
Since all $\alpha_j$ are positive and included in the $\mathcal O$ notation in \eqref{gen_J}, we simplify the notation by letting all $\alpha_j=1$ below.  The loss function
\[
\begin{split}
\Loss\big(\bar{v}(x),v(x)\big):=|\bar{v}(x)-v(x)|^2 +|\bar{v}'(x)-v'(x)|^2\,,
\end{split}
\]
satisfies
\begin{equation}\label{LL}
\begin{split}
&\Loss\big(\bar{v}(x),v(x)\big)= 
\sum_{k=1}^K\sum_{\ell=1}^K\eta_k\eta_\ell^*e^{\mathrm{i}(\omega_k-\omega_\ell)\cdot x} + |v(x)|^2
\\
& \quad -\sum_{k=1}^K\big(\eta_ke^{\mathrm{i}\omega_k\cdot x}v(x)^*+\eta_k^*e^{-\mathrm{i}\omega_k\cdot x}
v(x)\big)+\sum_{k=1}^K\sum_{\ell=1}^K\eta_k\eta_\ell^*e^{\mathrm{i}(\omega_k-\omega_\ell)\cdot x}\omega_k\cdot\omega_\ell +|v'(x)|^2
\\
&  \quad -\sum_{k=1}^K\big(\eta_ke^{\mathrm{i}\omega_k\cdot x}\mathrm{i}\omega_k\cdot (v')^*(x)
-\eta_k^*e^{-\mathrm{i}\omega_k\cdot x}\mathrm{i}\omega_k\cdot v'(x)\big)\,.\\
\end{split}
\end{equation}
Following the notations defined in \eqref{def_Ex_EJ}
\[
\begin{split}
\mathbb E_x[\Loss\big(\bar{v}(x),v(x)\big)]&=\int_{\rset^{3n}}\Loss\big(\bar{v}(x),v(x)\big)\mu_x(x)\mathrm{d}x\,,\\
\widehat{\mathbb E}_J[\Loss\big(\bar{v}(x),v(x)\big)]&=\frac{1}{J}\sum_{j=1}^J\Loss
\big(\bar{v}(x_j),v(x_j)\big)\,,\\
\end{split}
\]
the regularization terms will provide bounds for $\|\eta\|_1$, $\|\eta\|^2_2$ and  $\|\eta\|^4_2$, so that we can estimate  the data approximation error
$|\widehat{\mathbb E}_J[e^{\mathrm{i}(\omega_k-\omega_\ell)\cdot x}]- \mathbb E_x[e^{\mathrm{i}(\omega_k-\omega_\ell)\cdot x}]|$ in \eqref{LL}.

{\it Step 1} %
Let
\begin{equation}\label{Fourier}
\begin{split}
h(x)&=\int_{\mathbb R^{3n}} \widehat h(\omega)e^{\mathrm{i} \omega\cdot x}\mathrm{d}\omega\,,\ x\in\tset^{3n}\,,\\
\widehat h(\omega) &= \frac{1}{(2\pi)^{3n}}\int_{\tset^{3n}} h(x) e^{-\mathrm{i}\omega\cdot x}\mathrm{d}x\,,\ \omega\in\mathbb R^{3n}\,,
\end{split}
\end{equation}
denote the Fourier transform relation of $h:\rset^{3n}\to \mathbb C$ and $\widehat h:\rset^{3n}\to \mathbb C$.
The amplitude coefficients given by \eqref{eta_F}
\[
\hat\eta_k = \frac{\widehat{v}(\omega_k)}{K\PN(\omega_k)}\,,\ k=1,\ldots, K\,,
\]
define the minimal variance approximation
\[
\Bw := \sum_{k=1}^K \hat\eta_k\, e^{\mathrm{i}\omega_k\cdot x}\,,
\]
that satisfies the orthogonality relation
\begin{equation}\label{zero_mean}
\mathbb E_\omega[\Bw(x)-v(x)]=0
\end{equation}
since
\begin{equation}\label{mean_ok}
\begin{split}
\mathbb E_\omega[\Bw(x)]&= \mathbb E_\omega[\sum_{k=1}^K   
\frac{\widehat v(\omega_k)e^{\mathrm{i}\omega_k\cdot x}}{K\PN(\omega_k)}] =\int_{\mathbb R^{3n}}\widehat{v}(\omega)  e^{\mathrm{i}\omega\cdot x}\mathrm{d}\omega=v(x)\,.
\end{split}
\end{equation}

We have by \eqref{eta_F} and \eqref{mean_ok} %
 \begin{equation}\label{gen_error2}
 \begin{split}
\mathbb E_\omega[|\Bw(x)-v(x)|^2]
&= \frac{1}{K}\big( \int_{\mathbb R^{3n}}\frac{|\widehat v(\omega)|^2}{\PN(\omega)}\mathrm{d}\omega -|v(x)|^2\big)
=\mathcal O(K^{-1})\,,\\
\mathbb E_\omega[\Loss\big(\Bw(x),v(x)\big)] 
&= \frac{1}{K}\big( \int_{\mathbb R^{3n}}\frac{|\widehat v(\omega)|^2(1+|\omega|^2)}{\PN(\omega)}\mathrm{d}\omega -|v(x)|^2 -| v'(x)|^2\big)\,,
\end{split}
\end{equation}
since by the independence of $\omega_k$ and the unbiased property \eqref{zero_mean}
\[
\begin{split}
    \mathbb E_\omega [|\Bw(x)-v(x)|^2] 
    &=\frac{1}{K^2}\sum_{k=1}^K\sum_{\ell=1}^K
    \mathbb E_\omega [
    \big(\frac{\widehat v(\omega_k)}{\PN(\omega_k)}e^{\mathrm{i} \omega_k\cdot x}-v(x)\big)
    \big(\frac{\widehat v(\omega_\ell)}{\PN(\omega_\ell)}e^{\mathrm{i} \omega_\ell\cdot x}-v(x)\big)^*]\\
    &=\frac{1}{K^2}\sum_{k=1}^K
    \mathbb E_\omega [
    |\frac{\widehat v(\omega_k)}{\PN(\omega_k)}e^{\mathrm{i} \omega_k\cdot x}-v(x)|^2]\\
    &=\frac{1}{K}\big(\mathbb E_\omega [
    |\frac{\widehat v(\omega_k)}{\PN(\omega_k)}|^2]-|v(x)|^2\big)\\
    \end{split}
\]
and similarly for $\Loss(\Bw,v)$.

The variance of $\Loss(\Bw,v)$ is minimized by using the optimal probability density by \eqref{rho_star_optimal}
\[
\PN_\ast(\omega)=\frac{|\widehat v(\omega)|\sqrt{1+|\omega|^2}}{\|\widehat v(\omega)\sqrt{1+|\omega|^2}\|_{L^1}}
\]
and then by \eqref{gen_error2} it follows
\[
\mathbb E_\omega[\Loss\big(\Bw(x),v(x)\big)] =\frac{1}{K}( \|\widehat v(\omega)\sqrt{1+|\omega|^2}\|_{L^1}^2 -|v(x)|^2 -|v'(x)|^2)\,.
\]
We have assumed that the sampling density $\PN$ satisfies
\begin{equation}\label{pnb}
\sup_{\omega\in\rset^n}\frac{|\widehat v(\omega)|(1+|\omega|^2)}{\PN(\omega)}\le C\|\widehat v(1+|\omega|^2)\|_{L^1} \,,
\end{equation}
which implies
\begin{equation}\label{C''}
\|\Bw\|_{L^\infty}
\le\sum_{k=1}^K |\hat\eta_k|
=\sum_{k=1}^K\big|\frac{\widehat{v}(\omega_k)}{K\PN(\omega_k)}\big|
\le  C\|\widehat v(1+|\omega|^2)\|_{L^1} \,.
\end{equation}
We chose $\tilde{C}$ sufficiently large, so that $\tilde{C}>  C\|\widehat v(1+|\omega|^2)\|_{L^1}$ in \eqref{opt_J}.
 The optimization problem \eqref{opt_J}, %
  \eqref{C''}, \eqref{pnb} and \eqref{gen_error2} imply
that the expected training error satisfies %
\begin{equation}\label{tr}
\begin{split}
&\mathbb E_\omega\Big[\min_{\eta\in \mathbb C^K }\Big( %
\widehat{\mathbb E}_J[\Loss\big(\bar{v}(x),v(x)\big)]
 +\lambda_1\sum_{k=1}^K|\eta_k|^2 + \lambda_2(\sum_{k=1}^K|\eta_k|^2)^2\\
& \qquad\qquad\qquad\qquad \qquad\qquad \qquad\qquad\quad \ \quad+\lambda_3\max(\sum_{k=1}^K|\eta_k|-\tilde{C},0) 
\Big)\Big]\\
&\le \mathbb E_\omega\big[\widehat{\mathbb E}_J[\Loss\big(\bar{\bar v}(x),v(x)\big)] 
+\lambda_1\sum_{k=1}^K|\hat\eta_k|^2 + \lambda_2(\sum_{k=1}^K|\hat\eta_k|^2)^2
+\lambda_3\underbrace{\max(\sum_{k=1}^K|\hat\eta_k|-\tilde{C},0)}_{=0}\big]\\
&=\mathcal O(\frac{1}{K}+\frac{\lambda_1}{K}+ \frac{\lambda_2}{K^2})
\,,\\
\end{split}
\end{equation}
 and we obtain
\begin{equation}\label{reg_est}
\mathbb E_\omega\big[
\lambda_3\max(\sum_{k=1}^K|\eta_k|-\tilde{C},0) 
=\mathcal O(\frac{1}{K}+\frac{\lambda_1}{K}+ \frac{\lambda_2}{K^2})\,.\\
\end{equation}

The generalization error (i.e. the test error) can by \eqref{tr} be written
\begin{equation}\label{LLE}
\begin{split}
&\mathbb E_\omega\big[\mathbb E_x[\Loss\big(\bar{v}(x),v(x)\big)]\big]\\
&=\mathbb E_\omega\big[\widehat{\mathbb E}_J[\Loss\big(\bar{v}(x),v(x)\big)]\big]
+\mathbb E_\omega\big[\mathbb E_x[\Loss\big(\bar{v}(x),v(x)\big)]\big]
-\mathbb E_\omega\big[\widehat{\mathbb E}_J[\Loss\big(\bar{v}(x),v(x)\big)]\big]\\
&=\mathcal O(\frac{1}{K}+\frac{\lambda_1}{K}+ \frac{\lambda_2}{K^2})
+\mathbb E_\omega\big[\mathbb E_x[\Loss\big(\bar{v}(x),v(x)\big)]\big]
-\mathbb E_\omega\big[\widehat{\mathbb E}_J[\Loss\big(\bar{v}(x),v(x)\big)]\big]\,.
\end{split}
\end{equation}
It remains to estimate  the right hand side using  \eqref{LL}.

{\it Step 2}
We will use that the order of integration and summation can be interchanged, that is
\begin{equation}\label{EE}
\begin{split}
\mathbb E_\omega\big[\mathbb E_x[\Loss(\bar{v},v)]\big] &= \mathbb E_x\big[\mathbb E_\omega[\Loss(\bar{v},v)]\big]\,,\\
\mathbb E_\omega\big[\widehat{\mathbb E}_J[\Loss(\bar{v},v)]\big] &= \widehat{\mathbb E}_J\big[\mathbb E_\omega[\Loss(\bar{v},v)]\big]\,.\\
\end{split}
\end{equation}

Splitting the loss function into two terms 
\[  \Loss\big(\bar{v}(x),v(x)\big) =  \Loss_{1}(\bar{v}(x),v(x)\big) +  \Loss_2(\bar{v}(x),v(x)\big)\,,\]
where
\begin{equation}\label{def_two_losses}
\begin{aligned}
    & \Loss_1(\bar{v}(x),v(x)\big):=|\bar{v}(x)-v(x)|^2\,,\\
    & \Loss_2(\bar{v}(x),v(x)\big):=|\bar{v}'(x)-v'(x)|^2\,,\\
\end{aligned}
\end{equation}
we have from \eqref{LLE} the expected generalization error
\begin{equation}\label{LLE_with_two_losses}
\begin{aligned}
    &\quad\ \mathbb{E}_{\{x_j \}}\Big[ \mathbb E_\omega\big[\mathbb E_x[\Loss(\bar{v}(x),v(x))\,\big|\,\{x_j \}]\big] \Big] \\
    & = \mathcal O(\frac{1}{K}+\frac{\lambda_1}{K}+ \frac{\lambda_2}{K^2})\\
&\quad +\mathbb{E}_{\{x_j \}}\Big[\mathbb E_\omega\big[\mathbb E_x[\Loss\big(\bar{v}(x),v(x)\big)\,\big|\,\{x_j \}]\big]\Big]
-\mathbb{E}_{\{x_j \}}\Big[\mathbb E_\omega\big[\widehat{\mathbb E}_J[\Loss\big(\bar{v}(x),v(x)\big)\,\big|\,\{x_j \}]\big]\Big]\\
&=\mathcal O(\frac{1}{K}+\frac{\lambda_1}{K}+ \frac{\lambda_2}{K^2})\\
& \quad +\mathbb{E}_{\{x_j \}}\Big[\mathbb E_\omega\big[\mathbb E_x[\Loss_1\big(\bar{v}(x),v(x)\big)\,\big|\,\{x_j \}]\big]\Big]-\mathbb{E}_{\{x_j \}}\Big[\mathbb E_\omega\big[\widehat{\mathbb E}_J[\Loss_1\big(\bar{v}(x),v(x)\big)\,\big|\,\{x_j \}]\big]\Big]\\
& %
\quad + \mathbb{E}_{\{x_j \}}\Big[\mathbb E_\omega\big[\mathbb E_x[\Loss_2\big(\bar{v}(x),v(x)\big)\,\big|\,\{x_j \}]\big]\Big]-\mathbb{E}_{\{x_j \}}\Big[\mathbb E_\omega\big[\widehat{\mathbb E}_J[\Loss_2\big(\bar{v}(x),v(x)\big)\,\big|\,\{x_j \}]\big]\Big]\,.
\end{aligned}
\end{equation}
Introducing the following notation for the trigonometric activation function
\begin{equation}\label{sigma_activation}
    \sigma(\omega_k\cdot x):= e^{\mathrm{i}\omega_k\cdot x}\,,
\end{equation}
we rewrite the terms involving $\Loss_1\big(\bar{v},v\big)$ and $\Loss_2\big(\bar{v},v\big)$ with the amplitude and frequency coefficients $\{\eta_k\}_{k=1}^K$ and $\{\omega_k\}_{k=1}^K$. Specifically for the terms with $\Loss_1\big(\bar{v},v\big)$ we have
\begin{equation} \label{Loss_1_rewrite}
\begin{aligned}
    & \quad\textup{ }\  \mathbb{E}_{\{x_j \}}\Big[\mathbb E_\omega\big[\mathbb E_x[\Loss_1\big(\bar{v}(x),v(x)\big)\,\big|\,\{x_j \}]\big]\Big]-\mathbb{E}_{\{x_j \}}\Big[\mathbb E_\omega\big[\widehat{\mathbb E}_J[\Loss_1\big(\bar{v}(x),v(x)\big)\,\big|\,\{x_j \}]\big]\Big]\\
    & = \mathbb{E}_{\{x_j \}}\Big[ \mathbb{E}_\omega\big[ (\mathbb{E}_x - \widehat{\mathbb{E}}_{J})[\Loss_1\big(\bar{v}(x),v(x)\,\big|\,\{x_j \}] \big]  \Big]\\
    & = \mathbb{E}_{\{x_j \}}\Big[ \mathbb{E}_\omega\big[ (\mathbb{E}_x - \widehat{\mathbb{E}}_{J})[ \sum_{k=1}^K \sum_{\ell=1}^K \eta_k^\ast\,\eta_\ell \,\sigma^\ast(\omega_k\cdot x)\,\sigma(\omega_\ell\cdot x)]\big]\Big] \\
    & \quad -2 \mathbb{E}_{\{x_j \}}\Big[ \mathbb{E}_\omega\big[ (\mathbb{E}_x - \widehat{\mathbb{E}}_{J})[ \sum_{k=1}^K  \mathrm{Re}\big(\eta_k \, \sigma(\omega_k\cdot x)\big)v(x)]\big]\Big] \\
    &\quad + \mathbb{E}_{\{x_j \}}\Big[ \mathbb{E}_\omega\big[ (\mathbb{E}_x - \widehat{\mathbb{E}}_{J})[ |v(x)|^2]\big]\Big]\,.
\end{aligned}
\end{equation}
Applying Cauchy's inequality we obtain the estimate
\begin{equation} \label{Loss_1_estimate_term_1}
\begin{aligned}
      & \quad\ \textup{ } \big| \mathbb{E}_\omega\big[ \sum_{k=1}^K \sum_{\ell=1}^K \eta_k^\ast\, \eta_\ell\,\big( \widehat{\mathbb{E}}_J[ \sigma^\ast(\omega_k\cdot x)\sigma(\omega_\ell\cdot x) ] - \mathbb{E}_x[ \sigma^\ast(\omega_k\cdot x)\sigma(\omega_\ell\cdot x) ] \big) \big]    \big|\\
      & \leq \Big( \mathbb{E}_\omega\big[ \sum_{k=1}^K \sum_{\ell=1}^K |\eta_k\,\eta_\ell|^2 \big] \Big)^{\frac{1}{2}}\, \Big( \mathbb{E}_{\omega}\big[ \sum_{k=1}^K \sum_{\ell=1}^K \big| (\widehat{\mathbb{E}}_J - \mathbb{E}_x)[ \sigma^\ast(\omega_k\cdot x)\sigma(\omega_\ell\cdot x) ] \big|^2 \big] \Big)^{\frac{1}{2}}\\
      & \leq \frac{\lambda_2}{2} \mathbb{E}_\omega\big[ \big( \sum_{k=1}^K |\eta_k|^2 \big)^2 \big] + \frac{1}{2\lambda_2} \mathbb{E}_\omega\big[ 
\sum_{k=1}^K \sum_{\ell=1}^K \big| ( \widehat{\mathbb{E}}_J - \mathbb{E}_x )[ \sigma^\ast(\omega_k\cdot x)\sigma(\omega_\ell\cdot x) ] \big|^2 \big]\, ,
\end{aligned}
\end{equation}
where in the last inequality we use that for any two real numbers $a$ and $b$ with a positive parameter $\gamma$ it holds that $2|ab|\leq \frac{a^2}{\gamma}+\gamma b^2$. By \eqref{tr} we have for the first regularization penalty term
\begin{equation} \label{sum_squared_amplitudes_estimate}
    \frac{\lambda_2}{2} \mathbb{E}_{\{ x_j \}}\Big[ \mathbb{E}_\omega\big[ \big( \sum_{k=1}^K |\eta_k|^2 \big)^2 \big]  \Big] = \mathcal{O}\big( \frac{1}{K} + \frac{\lambda_1}{K} + \frac{\lambda_2}{K^2} \big)\,.
\end{equation}
Noticing that the trigonometric activation function $\sigma(\omega\cdot x)$ has finite moments
\begin{equation} \label{sigma_moments_bounded}
    \mathbb{E}_\omega\big[ \mathbb{E}_x[ |\sigma(\omega\cdot x)|^2 +  |\sigma(\omega\cdot x)|^4 ] \big]=\mathcal{O}(1)\,,
\end{equation}
we use the independence of the data $\{x_j \}_{j=1}^J$ and the frequencies $\{\omega_k\}_{k=1}^K$ together with the zero expected value
\[
\mathbb{E}_{\{x_j \}}[ \sigma^\ast(\omega_k\cdot x_j)\sigma(\omega_\ell\cdot x_j) ]-\mathbb{E}_{x}[\sigma^\ast(\omega_k\cdot x)\sigma(\omega_\ell\cdot x)]=0
\]
to estimate the last term in \eqref{Loss_1_estimate_term_1}. Specifically we obtain
\begin{equation}\label{estimate_Loss_1_squared_sigma}
\begin{aligned}
    & \quad \ \mathbb{E}_{\{x_j\}}\Big[ \mathbb{E}_\omega\big[ \sum_{k=1}^K \sum_{\ell=1}^K \big| (\widehat{\mathbb{E}}_{J}-\mathbb{E}_x)[\sigma^\ast(\omega_k\cdot x)\sigma(\omega_\ell\cdot x)] \big|^2 \big] \Big]\\
    & = \sum_{k=1}^K \sum_{\ell=1}^K \mathbb{E}_\omega\Big[ \mathbb{E}_{\{x_j\}}\big[ \big| \widehat{\mathbb{E}}_{J}[\sigma^\ast(\omega_k\cdot x)\sigma(\omega_\ell\cdot x)] -\underbrace{\mathbb{E}_x[\sigma^\ast(\omega_k\cdot x)\sigma(\omega_\ell\cdot x)]}_{=:E_{k,\ell}}  \big|^2 \big] \Big]\\
    & = \sum_{k=1}^K \sum_{\ell=1}^K \mathbb{E}_\omega\Big[ \mathbb{E}_{\{x_j\}}\big[ \sum_{j=1}^J \sum_{i=1}^J \frac{\big( \sigma^\ast(\omega_k\cdot x_j)\sigma(\omega_\ell\cdot x_j) - E_{k,\ell} \big)^\ast}{J}\times\\
    &\qquad \qquad \qquad \qquad \qquad \qquad \times \frac{\big( \sigma^\ast(\omega_k\cdot x_i)\sigma(\omega_\ell\cdot x_i) - E_{k,\ell} \big)}{J}   \big]\Big]\\
    & = \sum_{k=1}^K \sum_{\ell=1}^K \mathbb{E}_\omega\Big[ \mathbb{E}_{\{x_j\}}\big[ \frac{1}{J^2}\sum_{j=1}^J \big| \sigma^\ast(\omega_k\cdot x_j)\sigma(\omega_\ell\cdot x_j) - E_{k,\ell} \big|^2 \big]\Big]\\
    & = \frac{1}{J} \sum_{k=1}^K \sum_{\ell=1}^K \mathbb{E}_\omega\Big[ \mathbb{E}_{\{x_j\}}\big[ \big| \sigma^\ast(\omega_k\cdot x_j)\sigma(\omega_\ell\cdot x_j) - E_{k,\ell} \big|^2 \big]\Big] \\
    & = \frac{1}{J} \sum_{k=1}^K \sum_{\ell=1}^K \mathbb{E}_\omega\Big[ \mathbb{E}_{x}\big[ \big| \sigma^\ast(\omega_k\cdot x)\sigma(\omega_\ell\cdot x) - E_{k,\ell} \big|^2 \big]\Big] \\
    & = \mathcal{O}\big( \frac{K^2}{J} \big)\,,
\end{aligned}
\end{equation}
where in the last equality we use the property of bounded moments of $\sigma(\omega\cdot x)$ in \eqref{sigma_moments_bounded}. The estimate \eqref{Loss_1_estimate_term_1}  together with 
 \eqref{sum_squared_amplitudes_estimate} and \eqref{estimate_Loss_1_squared_sigma} implies
\begin{equation} \label{Loss_1_term_1_final_estimate}
\begin{aligned}
    & \quad\ \textup{ } \Big| \mathbb{E}_{\{x_j\}} \Big[\mathbb{E}_\omega\big[ \sum_{k=1}^K \sum_{\ell=1}^K \eta_k^\ast\, \eta_\ell\big( \widehat{\mathbb{E}}_J[ \sigma^\ast(\omega_k\cdot x)\sigma(\omega_\ell\cdot x) ] - \mathbb{E}_x[ \sigma^\ast(\omega_k\cdot x)\sigma(\omega_\ell\cdot x) ] \big) \big]    \Big]\Big|\\
    & = \mathcal{O}\big( \frac{1}{K}+\frac{\lambda_1}{K} +\frac{\lambda_2}{K^2}+\frac{K^2}{\lambda_2 J} \big)\,.
\end{aligned}
\end{equation}

The second term on the right hand side of \eqref{Loss_1_rewrite} can be estimated similarly with Cauchy's inequality as
\begin{equation} \label{Loss_1_estimate_term_2}
\begin{aligned}
     & \quad\ \textup{ } \Big|\mathbb{E}_{\{x_j\}} \Big[  \mathbb{E}_\omega\big[ (\widehat{\mathbb{E}}_J - \mathbb{E}_x)[ \sum_{k=1}^K \eta_k\, \sigma(\omega_k\cdot x) \,v(x) ]  \big]  \Big] \Big|\\
     & = \Big|\mathbb{E}_{\{x_j\}} \Big[  \mathbb{E}_\omega\big[ \sum_{k=1}^K \eta_k\, (\widehat{\mathbb{E}}_J - \mathbb{E}_x)[ \sigma(\omega_k\cdot x) \,v(x) ] \big] \Big] \Big|\\
     & \le  \mathbb{E}_{\{x_j\}} \Big[ \Big( \mathbb{E}_\omega \big[ \sum_{k=1}^K |\eta_k|^2 \big] \Big)^{\frac{1}{2}}\, \Big( \mathbb{E}_\omega\big[ \sum_{k=1}^K \big| (\widehat{\mathbb{E}}_J - \mathbb{E}_x)[\sigma(\omega_k\cdot x)\,v(x)] \big|^2 \big] \Big)^{\frac{1}{2}} \Big]\\
     & \leq \frac{\lambda_1}{2}\mathbb{E}_{\{x_j\}} \Big[ \mathbb{E}_\omega \big[ \sum_{k=1}^K |\eta_k|^2 \big] \Big] + \frac{1}{2\lambda_1} \mathbb{E}_{\{x_j\}} \Big[ \mathbb{E}_\omega \big[ \sum_{k=1}^K \big| (\widehat{\mathbb{E}}_J - \mathbb{E}_x)[ \sigma(\omega_k\cdot x)\, v(x)] \big|^2  \big] \Big]\,,
\end{aligned}
\end{equation}
where for the first regularization penalty term we have by \eqref{tr} the estimate
\begin{equation}\label{sum_squared_amplitudes_estimate_lambda_1}
    \frac{\lambda_1}{2} \mathbb{E}_{\{ x_j \}}\Big[ \mathbb{E}_\omega\big[ \sum_{k=1}^K |\eta_k|^2 \big]  \Big] = \mathcal{O}\big( \frac{1}{K} + \frac{\lambda_1}{K} + \frac{\lambda_2}{K^2} \big)\,.
\end{equation}
As in \eqref{estimate_Loss_1_squared_sigma}, the independence of $\{ x_j \}_{j=1}^J$ and $\{\omega_k\}_{k=1}^K$ implies that
\begin{equation} \label{estimate_Loss_1_sigma_times_v1}
\begin{aligned}
    & \quad \textup{ }\ \frac{1}{2\lambda_1} \mathbb{E}_{\{x_j\}} \Big[ \mathbb{E}_\omega \big[ \sum_{k=1}^K \big| (\widehat{\mathbb{E}}_J - \mathbb{E}_x)[ \sigma(\omega_k\cdot x)\, v(x)] \big|^2  \big] \Big]\\
    & = \frac{1}{2\lambda_1} \mathbb{E}_\omega \Big[ \mathbb{E}_{\{x_j\}} \big[ \sum_{k=1}^K \big| (\widehat{\mathbb{E}}_J - \mathbb{E}_x)[ \sigma(\omega_k\cdot x)\, v(x)] \big|^2  \big] \Big]\\
    & = \frac{1}{2\lambda_1} \mathbb{E}_\omega \Big[ \sum_{k=1}^K \mathbb{E}_{\{x_j\}}\big[ \frac{1}{J^2}\sum_{j=1}^J \sum_{i=1}^J \big( \sigma(\omega_k\cdot x_j )\, v(x_j)-\mathbb{E}_x[\sigma(\omega_k\cdot x )\, v(x)] \big)^\ast\times \\
    & \phantom{\frac{1}{2\lambda_1} \mathbb{E}_\omega \Big[ \sum_{k=1}^K \mathbb{E}_{\{x_j\}}\big[ \frac{1}{J^2}\sum_{j=1}^J \sum_{i=1}^J}\qquad \times \big( \sigma(\omega_k\cdot x_i )\, v(x_i)-\mathbb{E}_x[\sigma(\omega_k\cdot x )\, v(x)] \big) \big]   \Big]\\
    & = \frac{1}{2\lambda_1} \mathbb{E}_\omega \Big[ \sum_{k=1}^K \mathbb{E}_{\{x_j\}}\big[ \frac{1}{J^2}\sum_{j=1}^J \big| \sigma(\omega_k\cdot x_j )\, v(x_j)-\mathbb{E}_x[\sigma(\omega_k\cdot x )\, v(x)] \big|^2 \big] \Big]\\
    & = \frac{K}{2\lambda_1 J} \mathbb{E}_\omega \Big[ \mathbb{E}_{\{x_j\}}\big[ \big| \sigma(\omega\cdot x_j )\, v(x_j)-\mathbb{E}_x[\sigma(\omega\cdot x )\, v(x)] \big|^2 \big] \Big]\\
    & = \mathcal{O}\big( \frac{K}{\lambda_1 J} \big)\,,
\end{aligned}
\end{equation}
where we also use the property of bounded second moment of the activation function $\sigma(\omega\cdot x)$ in \eqref{sigma_moments_bounded}. The estimate \eqref{Loss_1_estimate_term_2} together with \eqref{sum_squared_amplitudes_estimate_lambda_1} and \eqref{estimate_Loss_1_sigma_times_v1} implies
\begin{equation} \label{Loss_1_term_2_final_estimate}
    \Big|\mathbb{E}_{\{x_j \}}\Big[ \mathbb{E}_\omega\big[\sum_{k=1}^K  \eta_k \big(\mathbb{E}_x[ \sigma(\omega_k\cdot x)v(x)] - \widehat{\mathbb{E}}_{J}[ \sigma(\omega_k\cdot x)v(x)]\big)\big]\Big]\Big|\\
    = \mathcal{O}\big( \frac{1}{K}+\frac{\lambda_1}{K}+\frac{K}{\lambda_1 J} +\frac{\lambda_2}{K^2} \big)\,.
\end{equation}
Also, we have the zero expected value
\begin{equation}\label{Loss_1_term_3_final_estimate}
    \mathbb{E}_{\{x_j\}}\big[ \mathbb{E}_x[|v(x)|^2] -\widehat{\mathbb{E}}_J[|v(x)|^2]  \big] = 0\,.
\end{equation}
The estimates \eqref{Loss_1_term_1_final_estimate}, \eqref{Loss_1_term_2_final_estimate}, and \eqref{Loss_1_term_3_final_estimate} together give the bound for \eqref{Loss_1_rewrite} as
\begin{equation} \label{Loss_1_final_estimate}
\begin{aligned}
    & \quad \textup{ } \ \mathbb{E}_{\{x_j \}}\Big[\mathbb E_\omega\big[\mathbb E_x[\Loss_1\big(\bar{v}(x),v(x)\big)\,\big|\,\{x_j \}]\big]\Big]-\mathbb{E}_{\{x_j \}}\Big[\mathbb E_\omega\big[\widehat{\mathbb E}_J[\Loss_1\big(\bar{v}(x),v(x)\big)\,\big|\,\{x_j \}]\big]\Big]\\
    & = \mathcal{O}\big( \frac{1}{K}+\frac{\lambda_1}{K} +\frac{K}{\lambda_1 J} +\frac{\lambda_2}{K^2} + \frac{K^2}{\lambda_2 J} \big)\,.\\
\end{aligned}
\end{equation}
Furthermore, we notice that the terms with $\Loss_2\big(\bar{v},v\big)$ in \eqref{LLE_with_two_losses} can be rewritten as
\begin{equation} \label{Loss_2_rewrite}
    \begin{aligned}
    & \quad\textup{ }\  \mathbb{E}_{\{x_j \}}\Big[\mathbb E_\omega\big[\mathbb E_x[\Loss_2\big(\bar{v}(x),v(x)\big)\,\big|\,\{x_j \}]\big]\Big]-\mathbb{E}_{\{x_j \}}\Big[\mathbb E_\omega\big[\widehat{\mathbb E}_J[\Loss_2\big(\bar{v}(x),v(x)\big)\,\big|\,\{x_j \}]\big]\Big]\\
    & = \mathbb{E}_{\{x_j \}}\Big[ \mathbb{E}_\omega\big[ (\mathbb{E}_x - \widehat{\mathbb{E}}_{J})[\Loss_2\big(\bar{v}(x),v(x)\big)\,\big|\,\{x_j \}] \big]  \Big]\\
    & =  \mathbb{E}_{\{x_j \}}\Big[ \mathbb{E}_\omega\big[ (\mathbb{E}_x - 
    \widehat{\mathbb{E}}_{J})[ |\sum_{k=1}^K \mathrm{i}\, \omega_k \,\eta_k\, e^{\mathrm{i}\omega_k\cdot x} - v'(x)|^2    ] \big] \Big]\,,\\
    \end{aligned}
\end{equation}
where we use the gradient of the random Fourier feature network 
$\Bar{v}'(x)=\sum_{k=1}^K \mathrm{i}\, \omega_k \, \eta_k\, e^{\mathrm{i}\omega_k\cdot x}$.

By introducing the activation function $\Tilde{\sigma}(\omega,x)$ with multi-dimensional output
\begin{equation}\label{sigma_tilde}
    \Tilde{\sigma}(\omega,x) := \mathrm{i}\,\omega_k\, e^{\mathrm{i}\omega_k\cdot x}\,,
\end{equation}
we understand the $\Loss_2\big(\bar{v},v\big)$ term in \eqref{Loss_2_rewrite} as the squared loss for the function $\sum_{k=1}^K \eta_k\,\Tilde{\sigma}(\omega,x)$ with amplitudes $\{\eta_k\}_{k=1}^K$. Under the assumption \eqref{pn}, 
the estimate for the terms of $\Loss_1\big(\bar{v},v\big)$ applies similarly to the terms of $\Loss_2\big(\bar{v},v\big)$ in \eqref{Loss_2_rewrite}, with the activation function $\Tilde{\sigma}(\omega,x)$. We arrive at the loss estimation
\begin{equation} \label{Loss_2_final_estimate}
\begin{aligned}
    & \quad \textup{ } \ \mathbb{E}_{\{x_j \}}\Big[\mathbb E_\omega\big[\mathbb E_x[\Loss_2\big(\bar{v}(x),v(x)\big)\,\big|\,\{x_j \}]\big]\Big]-\mathbb{E}_{\{x_j \}}\Big[\mathbb E_\omega\big[\widehat{\mathbb E}_J[\Loss_2\big(\bar{v}(x),v(x)\big)\,\big|\,\{x_j \}]\big]\Big]\\
    & = \mathcal{O}\big( \frac{1}{K}+\frac{\lambda_1}{K} +\frac{K}{\lambda_1 J} +\frac{\lambda_2}{K^2} + \frac{K^2}{\lambda_2 J} \big)\,.\\
\end{aligned}
\end{equation}

{\it Step 3}
By combining \eqref{LLE}, \eqref{LLE_with_two_losses}, \eqref{Loss_1_final_estimate}, and \eqref{Loss_2_final_estimate}, we obtain
\begin{equation}
    \mathbb{E}_{\{x_j \}}\Big[ \mathbb E_\omega\big[\mathbb E_x[\Loss(\bar{v}(x),v(x))\,\big|\,\{x_j \}]\big] \Big]=\mathcal{O}\big( \frac{1}{K}+\frac{\lambda_1}{K} +\frac{K}{\lambda_1 J} +\frac{\lambda_2}{K^2} + \frac{K^2}{\lambda_2 J} \big)\,,
\end{equation}
which provides the generalization error bound
\[
\begin{aligned}
&\quad \mathbb E_{\{x_j\}}\Big[\mathbb E_\omega\big[%
\int_{\rset^{3n}}\big(\alpha_1|v(x)-\bar{v}(x)|^2 + \alpha_2|v'(x)-\bar{v}'(x)|^2\big)\mu_x(x)\,\mathrm{d}x\big]\Big]\\
&=\mathcal{O}\big( \frac{1}{K}+\frac{\lambda_1}{K} +\frac{K}{\lambda_1 J} +\frac{\lambda_2}{K^2} + \frac{K^2}{\lambda_2 J} \big)\,.
\end{aligned}
\]
We notice that the choice of parameters $\lambda_2=CK^2J^{-1/2}$ and $\lambda_1=CKJ^{-1/2}$, for some positive constant $C$, yields
\[
    \mathbb E_{\{x_j\}}\Big[\mathbb E_\omega\big[%
\int_{\rset^{3n}}\big(\alpha_1|v(x)-\bar{v}(x)|^2 + \alpha_2|v'(x)-\bar{v}'(x)|^2\big)\mu_x(x)\,\mathrm{d}x\big]\Big]
=\mathcal O\big(\frac{1}{K}+\sqrt{\frac{1}{J}}\big)\,.
\]

\end{proof}
\begin{preremark}
We note that the crucial cancellation in \eqref{estimate_Loss_1_squared_sigma} for $i\ne j$ uses that the frequencies $\{\omega_k\}_{k=1}^K$ are independent of the data $\{x_j \}_{j=1}^J$ which holds for random feature approximations but not for neural networks based on minimization over both the amplitudes $\{\eta_k\}_{k=1}^K$ and the frequencies $\{ \omega_k \}_{k=1}^K$.
\end{preremark}

 \subsection{Proof of Theorem \ref{thm}}\label{proof_main_thm}
 \if\JOURNAL1 
 \leavevmode \\
 \fi
The proof has four steps: 
\begin{description}
 \item[\it Step 1] derives  an error representation for observables based on a perturbed potential, using the global error given as a weighted sum of local errors along the approximate path with weights obtained from a transport equation for the exact path; 
\item[\it Step 2] derives an error estimate for derivatives of observables using the flow;
\item[\it Step 3] estimates the observable by combining Steps 1 and 2;
\item[\it Step 4] combines the remainder terms in Step 3 to obtain the bound \eqref{C_b}.
\end{description}

\begin{proof}

{\it Step 1 (Error representation).} Based on the derivation in Section~\ref{global_local_error},
we will use \eqref{g_est} with $g=A$ in the case that 
$\sum_{k=1}^K|\eta_k|< \tilde{C}+1$,
which by \eqref{reg_est} holds with probability $1-\mathcal O\big(\frac{1}{K}+\frac{\lambda_1}{K}+ \frac{\lambda_2}{K^2}\big)$, using that $\lambda_3=1$.
In the other case $\sum_{k=1}^K|\eta_k|\ge \tilde{C}+1$,
that holds with probability $\mathcal O\big(\frac{1}{K}+\frac{\lambda_1}{K}+ \frac{\lambda_2}{K^2}\big)$ with respect to $\omega$,
we have $|C_{AB}-\bar C_{AB}|=\mathcal O(1)$ by assumption \eqref{pn1}.

The generalization error estimate of $v-\bar{v}$ in Theorem \ref{thm:generalization}  can be applied to bound $V-\bar V$ by using the localization property $\chi_1\chi_0=\chi_0$ to obtain
\begin{equation}\label{vvv}
\begin{split}
V-\bar V&= V-\big(\bar{v}\chi_0+V(1-\chi_0)\big)\\
&=V\chi_0-\bar{v}\chi_0\\
&=V\chi_1\chi_0-\bar{v}\chi_0\\
&=(v-\bar{v})\chi_0\,,
\end{split}
\end{equation}
in combination with the property that $\chi_0$ and its derivatives are bounded. We note that 
\[|V-\bar v_r|=|V-{\rm Re}(\bar V)|\le |V-\bar V|\,.\]

{\it Step 2 (Error representation of gradients).} To verify that $ u'$ is bounded
we use the flow
\[
\frac{\partial g( z_s)}{\partial  z_i(t)} = \sum_j \underbrace{\frac{\partial z_j(s)}{\partial z_i(t)}}_{=: z'_{ij}(s)} \frac{\partial g( z_s)}{\partial  z_j(s)}
\]
which satisfies
\[
\frac{\mathrm{d}}{\mathrm{d}s}  z'_{ij}(s) = \left\{
\begin{array}{cc} \sum_k  z'_{ik}(s) f'_{jk}( z_s) & s>t\\
0 & s<t\,,
\end{array}\right.\,, \ z'(0)=\rm{I}\,.
\]
We have
\[
\partial_{z_i} u(z,t)= \partial_{z_i(t)} g( z_T; z_t=z) 
= \sum_j \frac{\partial  z_j(T)}{\partial z_i(t)} 
\frac{\partial g( z_T; z_t=z)}{\partial  z_j(T)}\,,
\] 
%
and   obtain
\begin{equation}\label{zF}
\dot{{{\mathbf z}}}_s = F(\mathbf{ z}_s,s)\,,
\end{equation}
by defining
\[
\begin{split}
\mathbf{ z}&:= \left[\begin{array}{c} 
 z\\
 z'
\end{array}
\right]\,, \\ 
 F({\mathbf z},s) &:= \left[\begin{array}{c} 
 p\\
-V'(x)\\
 z' \left[
\begin{array}{cc}
0 & - V''( x)\\
\rm{I} & 0\end{array}
\right]
\end{array}
\right]\,, \ s>t\,,\\
F({\mathbf z},s) &:= \left[\begin{array}{c} 
 p\\
- V'( x)\\
0
\end{array}
\right]\,, \ s<t\,.\\
\end{split}
\]

The assumptions that $\| V'\|_{L^\infty}$ and $\| V''\|_{L^\infty}$
are bounded show
together with \eqref{zF}  that also
\begin{equation}\label{mub}
\|\partial_{z_i(t)}  u(z,t)\|_{L^\infty}
=\| z'(T)g'( z_T; z_t=z)\|_{L^\infty}=\mathcal O(1)
\end{equation} is bounded. 

{\it Step 3 (Estimates of the observable).} Consider the set
\[\mathcal V:=\{\omega : \sum_{k=1}^K|\eta_k|< \tilde{C}+1\}
\,\]
and its complement set $\mathcal V^c:=\{\omega : \sum_{k=1}^K|\eta_k|\ge \tilde{C}+1\}$.
Split the observable error
\[
\mathcal{C}_{AB} -\bar{\mathcal{C}}_{AB}=(\mathcal{C}_{AB} -\bar{\mathcal{C}}_{AB}) \mathbbm 1_{\mathcal V}
+(\mathcal{C}_{AB} -\bar{\mathcal{C}}_{AB}) \mathbbm 1_{\mathcal V^c}.
\]
Since $|\mathcal{C}_{AB} -\bar{\mathcal{C}}_{AB}|$ is uniformly bounded we have by 
Step 1 and $\lambda_1=KJ^{-1/2}$, $\lambda_2=K^2J^{-1/2}$ that 
\[
\mathbb E_{\{x_j\}}\big[\mathbb E_\omega[(\mathcal{C}_{AB} -\bar{\mathcal{C}}_{AB}) \mathbbm 1_{\mathcal V^c}]\big]=
\mathcal O\big(\frac{1}{K} 
+ (\frac{1}{J})^{1/2}\big)\,.
\]
It remains to estimate $\mathbb E_{\{x_j\}}\big[\mathbb E_\omega[(\mathcal{C}_{AB} -\bar{\mathcal{C}}_{AB}) \mathbbm 1_{\mathcal V}]\big]$. 
\begin{equation}\label{R12}
\begin{split}
(\mathcal{C}_{AB} -\bar{\mathcal{C}}_{AB}) \mathbbm 1_{\mathcal V} %
&=\int_{\rtset^{6n}} A\big(z_t(z_0)\big)B(z_0) \mu(z_0)\mathrm{d}z_0\mathbbm 1_{\mathcal V}
-\int_{\rtset^{6n}} A\big(\bar z_t(z_0)\big)B(z_0) \bar \mu(z_0)\mathrm{d}z_0\mathbbm 1_{\mathcal V}\\
&=\int_{\rtset^{6n}} \Big(A\big(z_t(z_0)\big)-A\big(\bar z_t(z_0)\big)\Big)B(z_0) \bar\mu(z_0)\mathrm{d}z_0\mathbbm 1_{\mathcal V}\\
&\quad+\int_{\rtset^{6n}} A\big( z_t(z_0)\big)B(z_0) \big(\mu(z_0)-\bar \mu(z_0)\big)\mathrm{d}z_0\mathbbm 1_{\mathcal V}\\
&=:R_1+R_2\,.
\end{split}
\end{equation}
The first term in the right hand side can be estimated using \eqref{error_of_g} with $g=A$ as follows
\[
\begin{split}
R_1&=\int_{\rtset^{6n}} \Big(A\big(z_t(z_0)\big)-A\big(\bar z_t(z_0)\big)\Big)B(z_0)\, \bar\mu(z_0)\,\mathrm{d}z_0
\mathbbm 1_{\mathcal V}\\
&=\int_{\rtset^{6n}} \int_0^t \big( \bar v_r'(\bar x_s)- V'(\bar x_s)\big)\cdot\nabla_p  u(\bar z_s,s) B(z_0)\,\bar\mu(z_0)\, \mathrm{d}s
\,\mathrm{d}z_0\mathbbm 1_{\mathcal V}\,.
\end{split}
\]
Gibbs density satisfies $\bar \mu(z_0)=\bar \mu(\bar z_s)$,  which  implies
\begin{equation}\label{R11}
\begin{split}
&\quad \; \big|\mathbb E_{\{x_j\}}\big[\mathbb E_\omega[R_1]\big]|\\
& =|\mathbb E_{\{x_j\}}\big[\mathbb E_\omega[\int_{\rtset^{6n}} \int_0^t \big(\underbrace{ \Bar{v}_r'(\bar x_s)-V'(\bar x_s)}_{=:\Delta v'(\bar x_s)}\big)\cdot\nabla_p  u(\bar z_s,s) B(z_0)
\,\underbrace{\bar \mu(z_0)}_{=\bar \mu(\bar z_s)}\, \mathrm{d}s
\,\mathrm{d}z_0\mathbbm 1_{\mathcal V}]\big]\big|\\
&=\big|\mathbb E_{\{x_j\}}\big[\mathbb E_\omega[\int_{\rtset^{6n}} \int_0^t \Delta v'(\bar x_s)\cdot\nabla_p  u(\bar z_s,s) B(z_0)
{\big(\mu(\bar z_s)+\underbrace{\Delta\mu(\bar z_s)}_{:=\bar \mu(\bar z_s)-\mu(\bar z_s)}\big)} \mathrm{d}s
\mathrm{d}z_0\mathbbm 1_{\mathcal V}]\big]\big|\\
&\le|\mathbb E_{\{x_j\}}\big[\mathbb E_\omega[\int_{\rtset^{6n}} \int_0^t \Delta v'(\bar x_s)\cdot\nabla_p  u(\bar z_s,s) B(z_0)
\mu(\bar z_s) \mathrm{d}s
\mathrm{d}z_0\mathbbm 1_{\mathcal V}]\big]|\\
&\quad+|\mathbb E_{\{x_j\}}\big[\mathbb E_\omega[\int_{\rtset^{6n}} \int_0^t \Delta v'(\bar x_s)\cdot\nabla_p  u(\bar z_s,s) B(z_0)
\Delta\mu(\bar z_s) \mathrm{d}s
\mathrm{d}z_0\mathbbm 1_{\mathcal V}]\big]|\\
&=:R_{11} + R_{12}\\
\end{split}
\end{equation}
where Theorem \ref{thm:generalization} together with \eqref{vvv} yields
\begin{equation}
\begin{split}
R_{11}&\le \big(\int_0^t\mathbb E_{\{x_j\}}\big[\mathbb E_\omega[ \int_{\rtset^{6n}} | \Delta v'(\bar x_s)|^2 \mu(\bar z_s)  \mathrm{d}z_0]\big]\mathrm{d}s\big)^{1/2}\\
&\qquad\times \big(\int_0^t \mathbb E_{\{x_j\}}\big[\mathbb E_\omega[\int_{\rtset^{6n}} | u'(\bar z_s,s)|^2 B^2(z_0) \mu(\bar z_s) \mathrm{d}z_0]\big]\mathrm{d}s\big)^{1/2}\\
&= \big(\int_0^t \mathbb E_{\{x_j\}}\big[\mathbb E_\omega[ \int_{\rtset^{6n}}| \Delta v'(\bar x_s)|^2{\mu(\bar z_s)} 
\underbrace{|\frac{\partial z_0}{\partial \bar z_s}|}_{=1} \mathrm{d}\bar z_s]\big]\mathrm{d}s\big)^{1/2}
\times\\
&\qquad\times
\big(\int_0^t \mathbb E_{\{x_j\}}\big[\mathbb E_\omega[ \int_{\rtset^{6n}} | \nabla_p u(\bar z_s,s)|^2 B^2(z_0)
\mu(\bar z_s)\underbrace{|\frac{\partial z_0}{\partial \bar z_s}|}_{=1} \mathrm{d}\bar z_s]\big] \mathrm{d}s\big)^{1/2}\\
&=\mathcal O\Big(\big(\frac{1}{K} 
+ (\frac{1}{J})^{1/2}\big)^{1/2}\Big)\,,
\end{split}
\end{equation}
where we use that $\|\nabla_z u\|_{L^\infty}=\mathcal O(1)$, as proved in \eqref{mub},
together with the assumption 
\if\JOURNAL1
\\
\fi
$\|B\|_{L^\infty}=\mathcal O(1)$.

To estimate $R_{2}$ and $R_{12}$ in \eqref{R12} and \eqref{R11}
we will use the estimate
\begin{equation}\label{mu_est}
\begin{split}
&\quad |\frac{e^{-\beta V(x)}}{\int_{\rset^{3n}}e^{-\beta V(x)} \mathrm{d}x} 
-\frac{e^{-\beta \bar v_r(x)}}{\int_{\rset^{3n}}e^{-\beta \bar v_r(x)} \mathrm{d}x}|\,\mathbbm 1_{\mathcal V}\\
&=\mathcal O(1)\beta |\Delta v(x)|\mu_x(x)
+\mathcal O(1)\mu_x(x)\int_{\rset^{3n}}\beta |\Delta v(x)|\mu_x(x)   \mathrm{d}x\,,\\
\end{split}
\end{equation}
which is proved in \eqref{mu_proof} below.

We have by Theorem \ref{thm:generalization},
\begin{equation}\label{R_12_bound}
\begin{split}
& \quad \;R_{12}=|\mathbb E_{\{x_j\}}\big[
\mathbb E_\omega[\int_{\rtset^{6n}} \int_0^t \Delta v'(\bar x_s)\cdot\nabla_p  u(\bar z_s,s) B(z_0)
\Delta\mu(\bar z_s) \mathrm{d}s
\underbrace{|\frac{\partial z_0}{\partial \bar z_s}|}_{=1} \mathrm{d}\bar z_s\mathbbm 1_{\mathcal V}]\big]|\\
&=\mathcal O(1)\mathbb E_{\{x_j\}}\big[
\mathbb E_\omega[\int_{\rset^{3n}}|\Delta v'(x)|\Big( |\Delta v(x)|\mu_x(x) 
+ \mu_x(x)\int_{\rset^{3n}} |\Delta v(x')|\mu_x(x') \mathrm{d}x'\Big)\mathrm{d}x]\big]\\
&=\mathcal O(1)\mathbb E_{\{x_j\}}\big[
\mathbb E_\omega[(\int_{\rset^{3n}} |\Delta v'(x)|^2\mu_x(x) \mathrm{d}x)^{1/2}(\int_{\rset^{3n}} |\Delta v(x)|^2\mu_x(x) \mathrm{d}x)^{1/2}]\big]\\
&=\mathcal O(1)\mathbb E_{\{x_j\}}\big[
\mathbb E_\omega[\int_{\rset^{3n}} |\Delta v'(x)|^2\mu_x(x) \mathrm{d}x 
+\int_{\rset^{3n}} |\Delta v(x)|^2\mu_x(x) \mathrm{d}x]\big]\\
&=\mathcal O\big(\frac{1}{K} 
+ (\frac{1}{J})^{1/2}\big)\,.
\end{split}
\end{equation}

The second term in \eqref{R12} has the bound
\[
\begin{split}
&|\mathbb E_{\{x_j\}}\big[\mathbb E_\omega[R_2]\big]|= | \mathbb E_{\{x_j\}}\big[\mathbb E_\omega[\int_{\rtset^{6n}} A\big(\bar z_t(z_0)\big)B(z_0) 
\big(\mu(z_0)-\bar \mu(z_0)\big)\mathrm{d}z_0\mathbbm 1_{\mathcal V}]\big]| \\
&\le \mathbb E_{\{x_j\}}\big[\mathbb E_\omega[\int_{\rtset^{6n}} |A\big(\bar z_t(z_0)\big)B(z_0)| 
\frac{e^{-\beta |p_0|^2/2}}{\int_{\rset^{3n}}e^{-\beta |p|^2/2} \mathrm{d}p}\times \\
& \qquad \qquad \qquad \quad \times \big|\frac{e^{-\beta V(x_0)}}{\int_{\rset^{3n}}e^{-\beta V(x)} \mathrm{d}x} 
-\frac{e^{-\beta \bar v_r(x_0)}}{\int_{\rset^{3n}}e^{-\beta \bar v_r(x)} \mathrm{d}x}\big|\,\mathrm{d}z_0\mathbbm 1_{\mathcal V}]\big]\\
&=\mathcal O(1)\|A\|_{L^\infty}\|B\|_{L^\infty}
\mathbb E_{\{x_j\}}\big[\mathbb E_\omega[\int_{\rtset^{6n}}|\frac{e^{-\beta V(x_0)}}{\int_{\rset^{3n}}e^{-\beta V(x)} \mathrm{d}x} 
-\frac{e^{-\beta \bar v_r(x_0)}}{\int_{\rset^{3n}}e^{-\beta \bar v_r(x)} \mathrm{d}x}|\mathrm{d}x_0\mathbbm 1_{\mathcal V}]\big]\,.\\
\end{split}
\]
We obtain by Theorem \ref{thm:generalization}, \eqref{vvv} and \eqref{mu_est}
\begin{equation}\label{R_22}
\begin{split}
|\mathbb E_{\{x_j\}}\big[\mathbb E_\omega[R_2]\big]|
&= \mathcal O(1)\mathbb E_{\{x_j\}}\big[\mathbb E_\omega[\int_{\rtset^{6n}} |\Delta v(x)|\mu_x(x)
\mathrm{d}x]\big]\\
&=\mathcal O(1)\big(\mathbb E_{\{x_j\}}\big[\mathbb E_\omega[\int_{\rtset^{6n}} |\Delta v(x)|^2\mu_x(x)\mathrm{d}x]\big]\big)^{1/2}\\
&=\mathcal O\Big(\big(\frac{1}{K}+(\frac{1}{J})^{1/2}\big)^{1/2}\Big)\,.
\end{split}
\end{equation}

{\it Proof of \eqref{mu_est}.} Theorem \ref{thm:generalization} and \eqref{vvv} yield
\[
\begin{split}
|e^{-\beta V(x)}-e^{-\beta \bar v_r(x)}|\mathbbm 1_{\mathcal V}&=e^{-\beta V(x)}|1-e^{-\beta\Delta v(x)}|\mathbbm 1_{\mathcal V}\\
&=e^{-\beta V(x)}\frac{|1-e^{-\beta\Delta v(x)}|}{\beta |\Delta v(x)|}\beta |\Delta v(x)|\mathbbm 1_{\mathcal V}\\
&\le \beta |\Delta v(x)|e^{-\beta V(x)} e^{\beta|\Delta v(x)|}\mathbbm 1_{\mathcal V}\\
&\le \beta |\Delta v(x)|e^{-\beta V(x)}e^{\beta\|\Delta v(x)\|_{L^\infty}}\mathbbm 1_{\mathcal V}\\
&\le \beta |\Delta v(x)|e^{-\beta V(x)} \mathcal O(1)\,,
\end{split}
\]
which implies
\begin{equation}\label{mu_proof}
\begin{split}
&|\frac{e^{-\beta V(x)}}{\int_{\rset^{3n}}e^{-\beta V(x)} \mathrm{d}x} 
-\frac{e^{-\beta \bar v_r(x)}}{\int_{\rset^{3n}}e^{-\beta \bar v_r(x)} \mathrm{d}x}|\mathbbm 1_{\mathcal V}\\
&= |\frac{e^{-\beta V(x)}}{\int_{\rset^{3n}}e^{-\beta V(x)} \mathrm{d}x} 
\big(1- e^{-\beta\Delta v(x)}(1+\frac{\int_{\rset^{3n}}e^{-\beta  V(x)}\mathrm{d}x}{\int_{\rset^{3n}}e^{-\beta \bar v_r(x)} \mathrm{d}x}-1)\big)|\mathbbm 1_{\mathcal V}\\
&=|\frac{e^{-\beta V(x)}}{\int_{\rset^{3n}}e^{-\beta V(x)} \mathrm{d}x}(1- e^{-\beta\Delta v(x)}) 
+\frac{e^{-\beta V(x)}e^{-\beta\Delta v(x)}}{\int_{\rset^{3n}}e^{-\beta V(x)} \mathrm{d}x} 
\frac{\int_{\rset^{3n}}(e^{-\beta  V(x)} - e^{-\beta \bar v_r(x)})\mathrm{d}x}{\int_{\rset^{3n}}e^{-\beta \bar v_r(x)} \mathrm{d}x}|\mathbbm 1_{\mathcal V}\\
&= \mathcal O(1)\beta |\Delta v(x)|e^{-\beta V(x)} 
+\mathcal O(1)e^{-\beta V(x)} \int_{\rset^{3n}}\beta |\Delta v(x)|e^{-\beta V(x)} \mathrm{d}x\\
&= \mathcal O(1)\beta |\Delta v(x)|\mu_x(x)
+\mathcal O(1)\mu_x(x) \int_{\rset^{3n}}\beta |\Delta v(x)|\mu_x(x) \mathrm{d}x\,.\\
\end{split}
\end{equation}

{\it Step 4 (Summation of remainder terms).}
The approximation error of the observable
has by Theorem \ref{thm:generalization}, \eqref{vvv}, \eqref{R11} and \eqref{R_22} the bound
\[
\begin{split}
&|\mathbb E_{\{x_j\}}\big[\mathbb E_\omega[\mathcal{C}_{AB}-\bar{\mathcal{C}}_{AB}]\big]|\\
&=|\mathbb E_{\{x_j\}}\big[\mathbb E_\omega[(\mathcal{C}_{AB}-\bar{\mathcal{C}}_{AB})\mathbbm 1_{\mathcal V^c}]\big]
+\mathbb E_{\{x_j\}}\big[\mathbb E_\omega[R_1+R_2]\big]|\\
&=\mathcal O\big((\frac{1}{K}+(\frac{1}{J})^{1/2})^{1/2}\big)\,.
\end{split}
\]
\end{proof}

\begin{preremark}
The potential can be rewritten as a function $\tilde V(r(x))=V(x)$ of all pair distances $r(x)=(r_{12},\ldots, r_{n-1 n})(x)$. This representation with respect to $r$ becomes automatically rotation and translation invariant.
It also seems advantageous with respect to localization, since for each particle a coordinate system based on
distances to particles in a neighbourhood is naturally formed.  Let
$R_{i\cdot}:=\frac{\partial r_\cdot}{\partial x_i}$, then the gradient $\nabla_r \tilde V$ can be determined
by the minimal norm solution $\nabla_r \tilde V= R^T(RR^T)^{-1}\nabla V$ of the underdetermined equation $\nabla_x \bar V(x)=R\nabla_r\bar{\tilde V}(r(x))$.
Consequently training  $\bar{\tilde V}(r)$ by approximations to $\tilde V(r)$ with respect to pair distanced seems promising using the loss function $\int_{\rset^{3n}}|\nabla_r \tilde V(r)-\nabla_r\bar{\tilde V}(r)|^2\mu_x(x)\mathrm{d}x$
and determine the forces using $\nabla_x \bar V(x)=R\nabla_r\bar{\tilde V}(r(x))$.
\end{preremark}

\begin{preremark}
To include the permutation symmetry of the potential when interchanging coordinates related to atoms of the same kind, it is useful to write the potential as the sum of the energies related to each atom $V=\sum_{i=1}^n V_i$, where $V_i=V_j$ if atom $i$ and atom $j$ are of the same kind \cite{behler, weinan}. The splitting of energy into energy per particle can be made by using perturbation of the electron eigenvalue problem when removing one particle as in \cite{pma}.
\end{preremark}

\begin{preremark}
The convergence rate in the theorems seems to hold
also for neural networks with other activation functions $\sigma:\rset\to \rset$, replacing $y\mapsto e^{\mathrm{i}y}$, provided \eqref{gen_error2} and \eqref{C''} hold
for the neural network $\Bw(x)=\sum_{k=1}^K\eta_k \sigma(\omega_k\cdot x+b_k)$, using the Barron representation \[v(x)=\int_{\rset^{3n+1}}\eta(\omega,b)\sigma\big(\omega\cdot x + b\big)\mathrm{d}\omega\mathrm{d}b\,,\]
with $\eta_k=\eta(\omega_k,b_k)\,, k=1,\ldots, K$. %
\end{preremark}

\section*{Data availability}
The \texttt{Matlab} code for the numerical implementations of Section~\ref{sec_numerics} is available through the following Github repository: \url{https://github.com/XinHuang2022/Random_Fourier_feature_NN_MD}. The repository contains the scripts and functions used for data collection, optimization, simulation, and visualization tasks performed in this work.

\section*{Acknowledgment}
This research was supported by
Swedish Research Council grant 2019-03725 and KAUST grant OSR-2019-CRG8-4033.3. 
The work of P.P. was supported in part by the U.S. Army Research Office Grant W911NF-19-1-0243. The computations were enabled by resources provided by the National Academic Infrastructure for Supercomputing in Sweden (NAISS) at PDC Center for High Performance Computing, KTH Royal Institute of Technology, partially funded by the Swedish Research Council through grant agreement no. 2022-06725.

\if\JOURNAL2
\bibliography{Bibliography_V1}
\fi

\if\JOURNAL1
\bibliographystyle{abbrv}

\fi

\end{document}